\documentstyle[leqno, amstex]{amsart}
\newtheorem{thm}{Theorem}[section]
\newtheorem{prop}[thm]{Proposition}
\newtheorem{Lemma}[thm]{Lemma}

\theoremstyle{definition}
\newtheorem{Problem}{Problem}[section]
\newtheorem{rem}[thm]{Remark}

\def\qed{~\rule{1mm}{2.5mm}}
\begin{document}

\title{Studies on the Garnier system in two variables \\}

\author{ Yusuke Sasano }
\maketitle

\begin{abstract}
We study some Hamiltonian structures of the Garnier system in two variables from the viewpoints of its symmetry and holomorphy properties. We also give a generalization of {\it Okamoto transformation \it}of the sixth Painlev\'e system.
\end{abstract}

\section{Introduction}
In this paper, we present a new symmetry of the Garnier system in two variables $t,s$ with the Hamiltonians $H_1,H_2$ (see \eqref{1} in Section 3) explicitly given by
\begin{equation}\label{extensionoka}
S:(q_1,p_1,q_2,p_2,t,s) \rightarrow \left(q_1+\frac{q_2p_2+\alpha_1}{p_1},p_1,\frac{p_2}{p_1},-q_2p_1,t,\frac{t}{s}\right).
\end{equation}
The transformation $S$ is birational and symplectic, but this system is not invariant under $S$.

By using \eqref{extensionoka}, we obtain a new expression of this system explicitly given as follows:
\begin{align}\label{uraS}
\begin{split}
dq_1&=\frac{\partial H'_1}{\partial p_1}dt+\frac{\partial H'_2}{\partial p_1}ds, \quad dp_1=-\frac{\partial H'_1}{\partial q_1}dt-\frac{\partial H'_2}{\partial q_1}ds,\\
dq_2&=\frac{\partial H'_1}{\partial p_2}dt+\frac{\partial H'_2}{\partial p_2}ds, \quad dp_2=-\frac{\partial H'_1}{\partial q_2}dt-\frac{\partial H'_2}{\partial q_2}ds,\\
H'_1 &=H_{VI}(q_1,p_1,t;\alpha_1+\alpha_4+\alpha_6,-\alpha_1-\alpha_2,\alpha_2,\alpha_1+\alpha_5,\alpha_1+\alpha_3)\\
&+\alpha_2\{\frac{(s-1)q_2}{(t-1)(t-s)}-\frac{sq_1}{t(t-s)}+\frac{q_1q_2}{t(t-1)}\}p_2+\alpha_4\frac{(q_1-q_2)p_1}{t-s}\\
&+\{\frac{2(s-1)q_1q_2}{(t-1)(t-s)}-\frac{tq_2^2+sq_1^2}{t(t-s)}+\frac{(q_1^2+t)q_2}{t(t-1)}\}p_1p_2,\\
H'_2 &=\pi(H'_1) \quad (2\alpha_1+\alpha_2+\dots+\alpha_6=1),
\end{split}
\end{align}
where the transformation $\pi$ is given by
\begin{align}
\begin{split}
\pi:(q_1,p_1,q_2,p_2,t,s;\alpha_1,\alpha_2, \dots ,\alpha_6) \rightarrow (q_2,p_2,q_1,p_1,s,t;\alpha_1,\alpha_4,\alpha_3,\alpha_2,\alpha_5,\alpha_6).
\end{split}
\end{align}
The symbol $H_{VI}$ denotes the Hamiltonian of $P_{VI}$ (see in Section 2). Comparing with the original ones, the form of each transformed Hamiltonian $H'_i \ (i=1,2)$ is $H_{VI}$ with new parameters in addition to new interaction term $R'_i \ (i=1,2)$.
\begin{rem} {\rm
On the polynomiality of the transformed Hamiltonian by \eqref{extensionoka}, the change of time-variables $t,s$ of this transformation is not essential. In fact, for the transformation without changing time-variables $t,s$
\begin{equation}
S':(q_1,p_1,q_2,p_2) \rightarrow \left(q_1+\frac{q_2p_2+\alpha_1}{p_1},p_1,\frac{p_2}{p_1},-q_2p_1\right),
\end{equation}
we can obtain a polynomial Hamiltonian.
\rm}
\end{rem}

\begin{figure}
\unitlength 0.1in
\begin{picture}(43.22,20.62)(14.39,-22.60)
%
\special{pn 20}%
\special{ar 3709 574 681 376  0.0000000 6.2831853}%
%
\special{pn 20}%
\special{ar 2120 1870 681 376  0.0000000 6.2831853}%
%
\special{pn 20}%
\special{ar 5080 1880 681 376  0.0000000 6.2831853}%
%
\special{pn 20}%
\special{pa 3220 840}%
\special{pa 2400 1520}%
\special{fp}%
%
\special{pn 20}%
\special{pa 4150 860}%
\special{pa 4730 1530}%
\special{fp}%
%
\special{pn 20}%
\special{pa 2790 1930}%
\special{pa 4380 1930}%
\special{fp}%
\put(33.3000,-5.5000){\makebox(0,0)[lb]{Invariant}}%
\put(33.6000,-7.2000){\makebox(0,0)[lb]{cycles}}%
\put(16.7000,-19.7000){\makebox(0,0)[lb]{Symmetries}}%
\put(46.4000,-18.6000){\makebox(0,0)[lb]{Holomorphy}}%
\put(46.5000,-20.3000){\makebox(0,0)[lb]{conditions}}%
%
\special{pn 20}%
\special{sh 0.600}%
\special{ar 2970 1420 8 22  0.0000000 6.2831853}%
%
\special{pn 20}%
\special{sh 0.600}%
\special{ar 2980 1660 8 22  0.0000000 6.2831853}%
\put(30.5000,-15.0000){\makebox(0,0)[lb]{Special solutions}}%
\put(30.5000,-17.4000){\makebox(0,0)[lb]{Bilinear forms}}%
\put(15.3000,-10.7000){\makebox(0,0)[lb]{Noumi-Yamada's work}}%
\put(16.1000,-12.5000){\makebox(0,0)[lb]{(see \cite{N2})}}%
\put(27.1000,-22.6000){\makebox(0,0)[lb]{Yamada-Sasano's work}}%
\put(27.2000,-24.3000){\makebox(0,0)[lb]{(see \cite{Sasa5})}}%
\put(43.7000,-10.6000){\makebox(0,0)[lb]{second center of}}%
\put(44.8000,-12.3000){\makebox(0,0)[lb]{blowing-ups}}%
\put(43.0000,-8.9000){\makebox(0,0)[lb]{Space of initial conditions:}}%
\put(46.0000,-13.9000){\makebox(0,0)[lb]{(see \cite{Sasa1,Sasa2,Sasa3})}}%
\end{picture}%
\label{fig:intro}
\caption{}
\end{figure}

As these results, for \eqref{uraS} we will obtain
\begin{enumerate}
\item Different type of holomorphy conditions from the Garnier system (see Section 5).

\item Different type of invariant cycles from the Garnier system (see Section 8).

\item New Dynkin diagram (see Figure 1).

\item New expression of birational symmetries with a natural extension of {\it Okamoto transformation} (see in Section 10).

\item Some algebraic solutions including an extension of {\it Umemura's solution} of $P_{VI}$ (see in Section 11).
\end{enumerate}
We also study the Garnier system in two variables from the viewpoint of its holomorphy and symmetry.

This paper is organized as follows. In Section 1, we review the sixth Painlev\'e system from the viewpoint of its holomorphy and symmetry. In section 2, we recall the Hamiltonians of the Garnier system in two variables and find the holomorphy conditions of this system. In Section 3, we will consider the relation between the holomorphy conditions with the accessible singularities of this system. In Section 4, we propose some new expressions of this system and generalizations of {\it Okamoto transformation \it}of the sixth Painlev\'e system. In Section 5, we will consider the relation between the holomorphy conditions with the accessible singularities of each system given in Section 4. In Section 6, we study some B{\"a}cklund transformations of the Garnier system in two variables. In Section 7, we study some B{\"a}cklund transformations of \eqref{uraS}. In the final section, we study some algebraic solutions for \eqref{uraS} including an extension of {\it Umemura's solution} of $P_{VI}$.

\section{Review of $P_{VI}$-case}
As is well-known, each of the Painlev\'e equations $P_J \ (J=VI,V,IV,III,II,I)$ is equivalent to a polynomial Hamiltonian system $H_J$. In 1997, K. Takano et al. studied some Hamiltonian structures of Painlev\'e systems except for the first Painlev\'e system (see \cite{T1}). They showed that each space of initial conditions of $H_J$ is obtained by gluing a finite number of copies of ${\Bbb C}^2$ via the birational and symplectic transformations. In each coordinate system, the Hamiltonian is expressed as a {\it polynomial} of the canonical coordinates. Moreover, they showed that by using these patching data in the phase space of each $H_J$, each Hamiltonian $H_J$ is uniquely determined.
For example, let us review the case of the sixth Painlev\'e system. The sixth Painlev\'e system can be written as the Hamiltonian system
\begin{align}\label{SPVI}
\begin{split}
&\frac{dq}{dt}=\frac{\partial H_{VI}}{\partial p}, \quad \frac{dp}{dt}=-\frac{\partial H_{VI}}{\partial q},\\
&H_{VI}(q,p,t;\alpha_0,\alpha_1,\alpha_2,\alpha_3,\alpha_4)\\
&=\frac{1}{t(t-1)}[p^2(q-t)(q-1)q-\{(\alpha_0-1)(q-1)q+\alpha_3(q-t)q\\
&+\alpha_4(q-t)(q-1)\}p+\alpha_2(\alpha_1+\alpha_2)(q-t)] \quad (\alpha_0+\alpha_1+2\alpha_2+\alpha_3+\alpha_4=1).
\end{split}
\end{align}
By the work of Okamoto, it is known that the system \eqref{SPVI} has the affine Weyl group symmetry of type $D_4^{(1)}$, whose generators $s_i, \ i=0,1,2,3,4$, are given by
\begin{align}\label{D4}
\begin{split}
s_0(q,p,t;\alpha_0,\alpha_1,\dots,\alpha_4) \rightarrow &(q,p-\frac{\alpha_0}{q-t},p,t;-\alpha_0,\alpha_1,\alpha_2+\alpha_0,\alpha_3,\alpha_4),\\
s_1(q,p,t;\alpha_0,\alpha_1,\dots,\alpha_4) \rightarrow &(q,p,t;\alpha_0,-\alpha_1,\alpha_2+\alpha_1,\alpha_3,\alpha_4),\\
s_2(q,p,t;\alpha_0,\alpha_1,\dots,\alpha_4) \rightarrow &(q+\frac{\alpha_2}{p},p,t;\alpha_0+\alpha_2,\alpha_1+\alpha_2,-\alpha_2,\alpha_3+\alpha_2,\alpha_4+\alpha_2),\\
s_3(q,p,t;\alpha_0,\alpha_1,\dots,\alpha_4) \rightarrow &(q,p-\frac{\alpha_3}{q-1},t;\alpha_0,\alpha_1,\alpha_2+\alpha_3,-\alpha_3,\alpha_4),\\
s_4(q,p,t;\alpha_0,\alpha_1,\dots,\alpha_4) \rightarrow &(q,p-\frac{\alpha_4}{q},t;\alpha_0,\alpha_1,\alpha_2+\alpha_4,\alpha_3,-\alpha_4).
\end{split}
\end{align}
The list \eqref{D4} should be read as
\begin{align*}
&s_0(q)=q, \quad s_0(p)=p-\frac{\alpha_0}{q-t}, \quad s_0(t)=t,\\
&s_0(\alpha_0)=-\alpha_0, \quad s_0(\alpha_1)=\alpha_1, \quad s_0(\alpha_2)=\alpha_2+\alpha_0,\\
&s_0(\alpha_3)=\alpha_3, \quad s_0(\alpha_4)=\alpha_4.
\end{align*}

The Hamiltonian $H_{VI}$ \eqref{SPVI} is a polynomial in the variables $q,p$. In this sense we call the system \eqref{SPVI} as a polynomial Hamiltonian system. Consider the following birational and symplectic transformations $r_i, \ i=0,1,2,3,4$:
\begin{align}\label{holoPVI}
\begin{split}
&r_0:x_0=-((q-t)p-\alpha_0)p,\ y_0=\frac{1}{p},\\
&r_1:x_1=\frac{1}{q}, \ y_1=-(qp+\alpha_1+\alpha_2)q, \\
&r_2:x_2=\frac{1}{q}, \ y_2=-(qp+\alpha_2)q,\\
&r_3:x_3=-((q-1)p-\alpha_3)p,\ y_3=\frac{1}{p},\\
&r_4:x_4=-(qp-\alpha_4)p,\ y_4=\frac{1}{p}.
\end{split}
\end{align}
Since the transformations $r_i$ are symplectic, the system \eqref{SPVI} is transformed into a Hamiltonian system, whose Hamiltonian may have poles. It is remarkable that the transformed system becomes again a polynomial Hamiltonian system for any $i=0,1,\dots,4$. Furthermore, this holomorphy property uniquely characterizes the system \eqref{SPVI}:

\begin{prop}
Let us consider a polynomial Hamiltonian system with Hamiltonian $H \in {\Bbb C}(t)[q,p]$. We assume that

$(A1)$ $deg(H)=5$ with respect to $q,p$.

$(A2)$ This system becomes again a polynomial Hamiltonian system in each coordinate system $r_i \ (j=0,1,\dots,4)$.

\noindent
Then such a system coincides with the system \eqref{SPVI}.
\end{prop}
In this paper, we call the conditions $r_i \ (j=0,1,..,4)$ {\it holomorphy conditions \it} of the sixth Painlev\'e system. The space of initial conditions of $P_{VI}$ is covered by these coordinate systems. Each coordinate system contains a one-parameter family of meromorphic solutions of this system.

\begin{rem}\rm{
If we look for a polynomial Hamiltonian system which admits the symmetry \eqref{D4}, we have to consider huge polynomial in variables $q,p,t,\alpha_i$. On the other hand, in the holomorphy requirement \eqref{holoPVI}, we only need to consider polynomials in $q,p$. This reduces the number of unknown coefficients drastically.
\rm}
\end{rem}

\section{Holomorphy condisions of the Garnier system in two variables}
Consider a Fuchsian differential equation on ${\Bbb P}^1$
\begin{equation}\label{Fuchsian}
\frac{d^2Y}{dZ^2}+P_1(Z)\frac{dY}{dZ}+P_2(Z)Y=0
\end{equation}
with regular singularities $Z=0,Z=1,Z=t,Z=s,Z=\infty$, apparent singularities $Z=q_1,Z=q_2$ and the Riemann scheme
\begin{equation}\label{scheme}
\begin{pmatrix}
Z=0 & Z=1 & Z=t & Z=s & Z=q_1 & Z=q_2 & Z=\infty\\
0 & 0 & 0 & 0 & 0 & 0 & \alpha_1\\
\alpha_6 & \alpha_5 & \alpha_3 & \alpha_4 & 2 & 2 & \alpha_1+\alpha_2
\end{pmatrix}
\end{equation}
assuming that the Fuchs relation
\begin{equation}
2\alpha_1+\alpha_2+\dots +\alpha_6=1
\end{equation}
is satisfied. The monodromy preserving deformations of the equation \eqref{Fuchsian} with the scheme \eqref{scheme} is described as a completely integrable Hamiltonian system. By using a transformation (see \cite{S}), we obtain the Garnier system in two variables, which is equivalent to the Hamiltonian system given by (see \cite{S}, cf. \cite{Oka})
\begin{align}\label{1}
\begin{split}
dq_1&=\frac{\partial H_1}{\partial p_1}dt+\frac{\partial H_2}{\partial p_1}ds, \quad dp_1=-\frac{\partial H_1}{\partial q_1}dt-\frac{\partial H_2}{\partial q_1}ds,\\
dq_2&=\frac{\partial H_1}{\partial p_2}dt+\frac{\partial H_2}{\partial p_2}ds, \quad dp_2=-\frac{\partial H_1}{\partial q_2}dt-\frac{\partial H_2}{\partial q_2}ds,\\
H_1 &=H_{VI}(q_1,p_1,t;\alpha_4+\alpha_6,\alpha_2,\alpha_1,\alpha_5,\alpha_3)\\
&+(2\alpha_1+\alpha_2)\frac{q_1q_2p_2}{t(t-1)}+\alpha_3\{\frac{p_1}{t-s}-\frac{(s-1)p_2}{(t-s)(t-1)}\}q_2+\alpha_4\frac{s(p_2-p_1)q_1}{t(t-s)}\\
&+\{\frac{2(s-1)p_1p_2}{(t-s)(t-1)}-\frac{tp_1^2+sp_2^2}{t(t-s)}+\frac{(2q_1p_1+q_2p_2)p_2}{t(t-1)}\}q_1q_2,\\
H_2&=\pi(H_1),
\end{split}
\end{align}
where the transformation $\pi$ is explicitly given by
\begin{align}
\begin{split}
\pi:&(q_1,p_1,q_2,p_2,t,s;\alpha_1,\alpha_2,\alpha_3,\alpha_4,\alpha_5,\alpha_6)\\
&\rightarrow(q_2,p_2,q_1,p_1,s,t;\alpha_1,\alpha_2,\alpha_4,\alpha_3,\alpha_5,\alpha_6).
\end{split}
\end{align}
The symbol $H_{VI}$ is given by \eqref{SPVI} in Section 0. We remark that Kimura and Okamoto introduced polynomial Hamiltonians for the Garnier system (see \cite{Oka}). After that, by addition to improve, finally it has been of the form $H_1,H_2$ by T. Tsuda (see \cite{Tsuda1,Tsuda2,Tsuda3,Tsuda4}).

Here we recall the definition of a symplectic transformation and its properties (see \cite{T1}). Let
$$
\varphi : x=x(X,Y,Z,W,t), \  y=y(X,Y,Z,W,t), \  z=z(X,Y,Z,W,t),
$$
$$
 \  w=w(X,Y,Z,W,t), \ t=t
$$
be a biholomorphic mapping from a domain $D$ in ${\Bbb C}^5 \ni (X,Y,Z,W,t)$ into ${\Bbb C}^5 \ni (x,y,z,w,t)$. We say that the mapping is symplectic if
\begin{equation*}
dx \wedge dy + dz \wedge dw = dX \wedge dY + dZ \wedge dW,
\end{equation*}
where $t$ is considered as a constant or a parameter, namely, if, for $t=t_0$, ${\varphi}_{t_0}=\varphi|_{t=t_0}$ is a symplectic mapping from the $t_0$-section $D_{t_0}$ of $D$ to $\varphi(D_{t_0})$. Suppose that the mapping is symplectic. Then any Hamiltonian system
$$
dx/dt=\partial H/\partial y, \ \  dy/dt=-\partial H/\partial x, \ \  dz/dt=\partial H/\partial w, \ \ dw/dt=-\partial H/\partial z
$$
is transformed to
$$
dX/dt=\partial K/\partial Y, \ \  dY/dt=-\partial K/\partial X, \ \ dZ/dt=\partial K/\partial W, \ \ dW/dt=-\partial K/\partial Z,
$$
where
\begin{equation*}\label{A}
(A) \ \ \ dx \wedge dy +dz \wedge dw - dH \wedge dt =dX \wedge dY +dZ \wedge dW - dK \wedge dt.
\end{equation*}
Here $t$ is considered as a variable. By this equation, the function $K$ is determined by $H$ uniquely modulo functions of $t$, namely, modulo functions independent of $X,Y,Z$ and $W$.

\begin{thm}\label{th;holoG}
Let us consider a polynomial Hamiltonian system with Hamiltonians $H_i \in {\Bbb C}(t,s)[q_1,p_1,q_2,p_2] \ (i=1,2)$. We assume that

$(A1)$ $deg(H_i)=5$ with respect to $q_1,p_1,q_2,p_2$.

$(A2)$ This system becomes again a polynomial Hamiltonian system in each coordinate $r_i, \ i=1,2,\dots,6${\rm:\rm}
\begin{align}
\begin{split}
&r_1:x_1=\frac{1}{q_1}, \ y_1=-q_1(q_1p_1+q_2p_2+\alpha_1), \ z_1=\frac{q_2}{q_1}, \ w_1=q_1p_2, \\
&r_2:x_2=\frac{1}{q_1}, \ y_2=-q_1(q_1p_1+q_2p_2+\alpha_1+\alpha_2), \ z_2=\frac{q_2}{q_1}, \ w_2=q_1p_2, \\
&r_3:x_3=-p_1(q_1p_1-\alpha_3),\ y_3=\frac{1}{p_1},\ z_3=q_2,\ w_3=p_2, \\
&r_4:x_4=q_1, \ y_4=p_1, \ z_4=-p_2(q_2p_2-\alpha_4), \ w_4=\frac{1}{p_2}, \\
&r_5:x_5=-((q_1+q_2-1)p_1-\alpha_5)p_1, \ y_5=\frac{1}{p_1}, \  z_5=q_2, \ w_5=p_2-p_1, \\
&r_6:x_6=-((q_1+tq_2/s-t)p_1-\alpha_6)p_1,\ y_6=\frac{1}{p_1}, \ z_6=q_2, \  w_6=p_2-\frac{tp_1}{s}.
\end{split}
\end{align}
Then such a system coincides with the system \eqref{1}.
\end{thm}
We remark that each transformation of each coordinate $r_i, \ i=1,2,\dots,6,$ is birational and symplectic. These transformations are appeared as the patching data in the space of initial conditions of the system \eqref{1}. On the construction of these transformations $r_i \ (i=1,2,\dots,6),$ we will explain in the next section.

\begin{prop} 
In each coordinate $r_i, \ i=1,2,\dots,6$, the  Hamiltonians $H_{j1}$ and $H_{j2}$ on $U_j \times B$ are expressed as a polynomial in $x_j,y_j,z_j,w_j$ and a rational function in $t$ and $s$, and satisfy the following conditions{\rm: \rm}
\begin{align}\label{symplectic}
\begin{split}
&dq_1 \wedge dp_1 +dz \wedge dp_2 - dH_1 \wedge dt- dH_2 \wedge ds\\
&=dx_j \wedge dy_j +dz_j \wedge dw_j - dH_{j1} \wedge dt- dH_{j2} \wedge ds \quad (j=1,2,\dots,5),\\
&dq_1 \wedge dp_1 +dq_2 \wedge dp_2 - d(H_1-(1-q_2/s)p_1) \wedge dt- d(H_2-(1-q_1/t)p_2) \wedge ds\\
&=dx_6 \wedge dy_6 +dz_6 \wedge dw_6 - dH_{61} \wedge dt- dH_{62} \wedge ds.
\end{split}
\end{align}
\end{prop}

{\bf Proof of Theorem \ref{th;holoG}.}
The polynomial $H$ satisfying $(A1)$ has 126 unknown coefficients in ${\Bbb C}(t,s)$.
At first, resolving the coordinate $r_1$ in the variables $q_1,p_1,q_2,p_2$, we obtain
\begin{equation}\tag*{($R_1$)}
q_1=1/x_1, \quad p_1=-(x_1y_1+z_1w_1+\alpha_1)x_1, \quad q_2=z_1/x_1, \quad p_2=w_1x_1.
\end{equation}
 By $R_1$, we transform $H$ into $R_1(H)$, which has poles in only $x_1$. For $R_1(H)$, we only have to determine the unknown coefficients so that they cancel the poles of $R_1(H)$.

For the transformations $r_i, (i=2,3,4,5)$, we can repeat the same way.

We must be careful of the case of the transformation $r_6$. The relation between the coordinate system $(x_6,y_6,z_6,w_6)$ and the coordinate system $(q_1,p_1,q_2,p_2)$ is given by
\begin{align}
\begin{split}
&dx_6 \wedge dy_6 +dz_6 \wedge dw_6\\
&=dq_1 \wedge dp_1 +dq_2 \wedge dp_2+d((1-q_2/s)p_1) \wedge dt+d((1-q_1/t)p_2) \wedge ds.
\end{split}
\end{align}
Resolving the coordinate $r_6$ in the variables $q_1,p_1,q_2,p_2$, we obtain
\begin{equation}\tag*{($R_6$)}
q_1=t-x_6y_6^2+\alpha_6 y_6-\frac{t z_6}{s}, \quad p_1=1/y_6, \quad q_2=z_6, \quad p_2=w_6+\frac{t}{sy_6}.
\end{equation}
In this case, we must consider the polynomiality for $R_6(H_1-(1-q_2/s)p_1)$ and $R_6(H_2-(1-q_1/t)p_2)$, respectively.

In this way, we can obtain the Hamiltonians $H_1,H_2$. \qed

\section{On some Hamiltonian structures of the system \eqref{1}}

In this section, we will consider the relation between the holomorphy conditions $r_j$  with the accessible singularities of the system \eqref{1}. Let us take the compactification
\begin{equation*}
(q_1,p_1,q_2,p_2,t,s) \in {\Bbb C}^4 \times B_2 \  {\rm{to\rm}} \ ([z_0:z_1:z_2:z_3:z_4],t,s) \in {\Bbb P}^4 \times B_2
\end{equation*}
with the natural embedding
\begin{equation*}
(q_1,p_1,q_2,p_2)=(z_1/z_0,z_2/z_0,z_3/z_0,z_4/z_0).
\end{equation*}
Here $B_2={\Bbb C}^2-\{t(t-1)s(s-1)=0 \}$. Fixing the parameters $\alpha_i$, consider the product ${\Bbb P}^4 \times B_2$ and extend the regular vector field on ${\Bbb C}^4 \times B_2$ to a rational vector field $\tilde{v}$ on ${\Bbb P}^4 \times B_2$. It is easy to see that ${\Bbb P}^4$ is covered by five copies of ${\Bbb C}^4${\rm : \rm}
$$
U_0={\Bbb C}^4 \ni (q_1,p_1,q_2,p_2), \ U_j={\Bbb C}^4 \ni (X_j,Y_j,Z_j,W_j) \ (j=1,2,3,4),
$$
via the following rational transformations
\begin{align}\label{coverofP4}
\begin{split}
&1)  X_1=1/q_1, \quad Y_1=p_1/q_1, \quad Z_1=q_2/q_1, \quad W_1=p_2/q_1,\\
&2)  X_2=q_1/q_2, \quad Y_2=p_1/q_2, \quad Z_2=1/q_2, \quad W_2=p_2/q_2,\\
&3)  X_3=q_1/p_1, \quad Y_3=1/p_1, \quad Z_3=q_2/p_1, \quad W_3=p_2/p_1,\\
&4)  X_4=q_1/p_2, \quad Y_4=p_1/p_2, \quad Z_4=q_2/p_2, \quad W_4=1/p_2.
\end{split}
\end{align}
By the following lemma, we will show that each coordinate system $(x_i,y_i,z_i,w_i)$ $(i=1,2,\dots,6)$ can be obtained by successive blowing-up procedures of the accessible singularities in the boundary divisor ${\mathcal H} (\cong {\Bbb P}^3) \subset {\Bbb P}^4$.

\begin{figure}
\unitlength 0.1in
\begin{picture}(37.59,18.30)(15.10,-23.80)
%
\special{pn 8}%
\special{pa 3490 550}%
\special{pa 1710 2320}%
\special{dt 0.045}%
\special{pa 1710 2320}%
\special{pa 1710 2320}%
\special{dt 0.045}%
%
\special{pn 8}%
\special{pa 1730 2330}%
\special{pa 5260 2340}%
\special{dt 0.045}%
\special{pa 5260 2340}%
\special{pa 5259 2340}%
\special{dt 0.045}%
%
\special{pn 8}%
\special{pa 3490 570}%
\special{pa 3490 1640}%
\special{dt 0.045}%
\special{pa 3490 1640}%
\special{pa 3490 1639}%
\special{dt 0.045}%
%
\special{pn 8}%
\special{pa 1730 2330}%
\special{pa 3490 1620}%
\special{dt 0.045}%
\special{pa 3490 1620}%
\special{pa 3489 1620}%
\special{dt 0.045}%
%
\special{pn 8}%
\special{pa 3490 1610}%
\special{pa 5250 2330}%
\special{dt 0.045}%
\special{pa 5250 2330}%
\special{pa 5249 2330}%
\special{dt 0.045}%
%
\special{pn 20}%
\special{pa 3480 560}%
\special{pa 3480 1630}%
\special{fp}%
%
\special{pn 20}%
\special{sh 0.600}%
\special{ar 1730 2320 29 19  0.0000000 6.2831853}%
%
\special{pn 20}%
\special{sh 0.600}%
\special{ar 2990 2320 29 19  0.0000000 6.2831853}%
%
\special{pn 20}%
\special{sh 0.600}%
\special{ar 5240 2340 29 19  0.0000000 6.2831853}%
%
\special{pn 20}%
\special{sh 0.600}%
\special{ar 3810 2330 29 19  0.0000000 6.2831853}%
%
\special{pn 8}%
\special{pa 3470 570}%
\special{pa 5250 2340}%
\special{dt 0.045}%
\special{pa 5250 2340}%
\special{pa 5250 2340}%
\special{dt 0.045}%
\put(35.0000,-11.5000){\makebox(0,0)[lb]{$C_1 \cup C_2$}}%
\put(15.1000,-25.3000){\makebox(0,0)[lb]{$P_3$}}%
\put(50.5000,-25.4000){\makebox(0,0)[lb]{$P_4$}}%
\put(27.7000,-25.4000){\makebox(0,0)[lb]{$P_5$}}%
\put(36.5000,-25.5000){\makebox(0,0)[lb]{$P_6$}}%
\end{picture}%
\label{fig:G1}
\caption{The figure denotes the boundary divisor ${\mathcal H}$ in ${\Bbb P}^4$. The dark parts correspond to the accessible singularities of the system \eqref{1}.}
\end{figure}
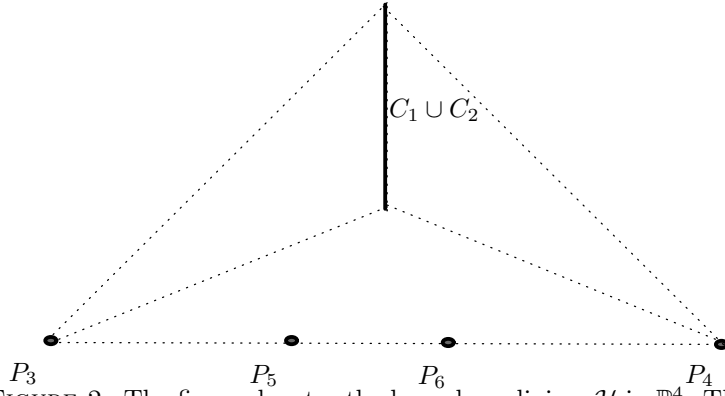

\begin{prop}\label{pro:asofG}
By resolving the following accessible singularities in the boundary divisor ${\mathcal H} (\cong {\Bbb P}^3) \subset {\Bbb P}^4$
\begin{equation*}
  \left\{
  \begin{aligned}
   C_1 &=\{(X_1,Y_1,Z_1,W_1)|X_1=Y_1=W_1=0 \},\\
   C_2 &=\{(X_2,Y_2,Z_2,W_2)|Y_2=Z_2=W_2=0 \},\\
   P_3 &=\{(X_3,Y_3,Z_3,W_3)|X_3=Y_3=Z_3=W_3=0\},\\
   P_4 &=\{(X_4,Y_4,Z_4,W_4)|X_4=Y_4=Z_4=W_4=0\},\\
   P_5 &=\{(X_3,Y_3,Z_3,W_3)|X_3=Y_3=Z_3=0,W_3=1\},\\
   P_6 &=\{(X_3,Y_3,Z_3,W_3)|X_3=Y_3=Z_3=0,W_3=t/s\},\\
   \end{aligned}
  \right. 
\end{equation*}
we can obtain the coordinate systems $(x_i,y_i,z_i,w_i)$ $(i=1,2,\dots,6)$. Here $C_1 \cup C_2 \cong {\Bbb P}^1$.
\end{prop}

\noindent
{\bf Proof of Proposition \ref{pro:asofG}.}
At first, we give an explicit resolution process for the accessible singularity $C_1 \cup C_2 \cong {\Bbb P}^1$ by giving a convenient coordinate system at each step.

By the following steps, we can resolve the accessible singularity $C_1 \cup C_2 \cong {\Bbb P}^1$.

{\bf Step 1:} We blow up along the curve $C_1 \cup C_2 \cong {\Bbb P}^1${\rm : \rm}
\begin{align*}
&{x_{1}^{(1)}}=X_1 \;, \;\;\; {y_{1}^{(1)}}=\frac{Y_1}{X_1} \;, \;\;\; {z_{1}^{(1)}}=Z_1 \;, \;\;\; {w_{1}^{(1)}}=\frac{W_1}{X_1},\\
&{x_{2}^{(1)}}=X_2 \;, \;\;\; {y_{2}^{(1)}}=\frac{Y_2}{Z_2} \;, \;\;\; {z_{2}^{(1)}}=Z_2 \;, \;\;\; {w_{2}^{(1)}}=\frac{W_2}{Z_2}.
\end{align*}

{\bf Step 2:} We blow up along the curve
\begin{align*}
&\{(x_{1}^{(1)},y_{1}^{(1)},z_{1}^{(1)},w_{1}^{(1)})|x_{1}^{(1)}=y_{1}^{(1)}=w_{1}^{(1)}=0\} \cup\\
&\{(x_{2}^{(1)},y_{2}^{(1)},z_{2}^{(1)},w_{2}^{(1)})|y_{2}^{(1)}=z_{2}^{(1)}=w_{2}^{(1)}=0\} \cong {\Bbb P}^1{\rm : \rm}
\end{align*}
\begin{align*}
&{x_{1}^{(2)}}=x_{1}^{(1)} \;, \;\;\; {y_{1}^{(2)}}=\frac{y_{1}^{(1)}}{x_{1}^{(1)}} \;, \;\;\; {z_{1}^{(2)}}=z_{1}^{(1)} \;, \;\;\; {w_{1}^{(2)}}=\frac{w_{1}^{(1)}}{x_{1}^{(1)}},\\
&{x_{2}^{(2)}}=x_{2}^{(1)} \;, \;\;\; {y_{2}^{(2)}}=\frac{y_{2}^{(1)}}{z_{2}^{(1)}} \;, \;\;\; {z_{2}^{(2)}}=z_{2}^{(1)} \;, \;\;\; {w_{2}^{(2)}}=\frac{w_{2}^{(1)}}{z_{2}^{(1)}}.
\end{align*}

It is easy to see that there are two accessible singularities
\begin{align*}
\begin{split}
S_1=&\{(x_{1}^{(2)},y_{1}^{(2)},z_{1}^{(2)},w_{1}^{(2)})|x_{1}^{(2)}=y_{1}^{(2)}+z_{1}^{(2)}w_{1}^{(2)}+\alpha_1=0\} \cup\\
&\{(x_{2}^{(2)},y_{2}^{(2)},z_{2}^{(2)},w_{2}^{(2)})|z_{2}^{(2)}=w_{2}^{(2)}+x_{2}^{(2)}y_{2}^{(2)}+\alpha_1=0\},\\
S_2=&\{(x_{1}^{(2)},y_{1}^{(2)},z_{1}^{(2)},w_{1}^{(2)})|x_{1}^{(2)}=y_{1}^{(2)}+z_{1}^{(2)}w_{1}^{(2)}+\alpha_1+\alpha_2=0\} \cup \\
&\{(x_{2}^{(2)},y_{2}^{(2)},z_{2}^{(2)},w_{2}^{(2)})|z_{2}^{(2)}=w_{2}^{(2)}+x_{2}^{(2)}y_{2}^{(2)}+\alpha_1+\alpha_2=0\}.
\end{split}
\end{align*}

{\bf Step 3:} We blow up along the surface
\begin{align*}
S_1=&\{(x_{1}^{(2)},y_{1}^{(2)},z_{1}^{(2)},w_{1}^{(2)})|x_{1}^{(2)}=y_{1}^{(2)}+z_{1}^{(2)}w_{1}^{(2)}+\alpha_1=0\} \cup\\
&\{(x_{2}^{(2)},y_{2}^{(2)},z_{2}^{(2)},w_{2}^{(2)})|z_{2}^{(2)}=w_{2}^{(2)}+x_{2}^{(2)}y_{2}^{(2)}+\alpha_1=0\}{\rm : \rm}
\end{align*}
\begin{align*}
&{x_{1}^{(3)}}=x_{1}^{(2)} \;, \;\;\; {y_{1}^{(3)}}=\frac{y_{1}^{(2)}+z_{1}^{(2)}w_{1}^{(2)}+\alpha_1}{x_{1}^{(2)}} \;, \;\;\; {z_{1}^{(3)}}=z_{1}^{(2)} \;, \;\;\; {w_{1}^{(3)}}=w_{1}^{(2)},\\
&{x_{2}^{(3)}}=x_{2}^{(2)} \;, \;\;\; {y_{2}^{(3)}}=y_{2}^{(2)} \;, \;\;\; {z_{2}^{(3)}}=z_{2}^{(2)} \;, \;\;\; {w_{2}^{(3)}}=\frac{w_{2}^{(2)}+x_{2}^{(2)}y_{2}^{(2)}+\alpha_1}{z_{2}^{(2)}}.
\end{align*}

{\bf Step 4:} We blow up along the surface
\begin{align*}
S_2=&\{(x_{1}^{(2)},y_{1}^{(2)},z_{1}^{(2)},w_{1}^{(2)})|x_{1}^{(2)}=y_{1}^{(2)}+z_{1}^{(2)}w_{1}^{(2)}+\alpha_1+\alpha_2=0\} \cup\\
&\{(x_{2}^{(2)},y_{2}^{(2)},z_{2}^{(2)},w_{2}^{(2)})|z_{2}^{(2)}=w_{2}^{(2)}+x_{2}^{(2)}y_{2}^{(2)}+\alpha_1+\alpha_2=0\}{\rm : \rm}
\end{align*}
\begin{align*}
&{x_{1}^{(4)}}=x_{1}^{(2)} \;, \;\;\; {y_{1}^{(4)}}=\frac{y_{1}^{(2)}+z_{1}^{(2)}w_{1}^{(2)}+\alpha_1+\alpha_2}{x_{1}^{(2)}} \;, \;\;\; {z_{1}^{(4)}}=z_{1}^{(2)} \;, \;\;\; {w_{1}^{(4)}}=w_{1}^{(2)},\\
&{x_{2}^{(4)}}=x_{2}^{(2)} \;, \;\;\; {y_{2}^{(4)}}=y_{2}^{(2)} \;, \;\;\; {z_{2}^{(4)}}=z_{2}^{(2)} \;, \;\;\; {w_{2}^{(4)}}=\frac{w_{2}^{(2)}+x_{2}^{(2)}y_{2}^{(2)}+\alpha_1+\alpha_2}{z_{2}^{(2)}}.
\end{align*}
We have resolved the accessible singularity $C_1 \cup C_2 \cong {\Bbb P}^1$.

By choosing new coordinate systems as
$$
(x_i,y_i,z_i,w_i)=(x_{1}^{(i+2)},-y_{1}^{(i+2)},z_{1}^{(i+2)},w_{1}^{(i+2)}) \ (i=1,2),
$$
we can obtain the coordinate systems $(x_i,y_i,z_i,w_i) \ (i=1,2)$ given in Theorem \ref{th;holoG}.

Next, we give an explicit resolution process for the accessible singular point $P_3$ by giving a convenient coordinate system at each step.

By the following steps, we can resolve the accessible singular point $P_3$.

{\bf Step 1:} We blow up at the point $P_3${\rm : \rm}
$$
{x_{3}^{(1)}}=\frac{X_3}{Y_3} \;, \;\;\; {y_{3}^{(1)}}=Y_3 \;, \;\;\; {z_{3}^{(1)}}=\frac{Z_3}{Y_3} \;, \;\;\; {w_{3}^{(1)}}=\frac{W_3}{Y_3}.
$$

{\bf Step 2:} We blow up along the surface
\begin{equation*}
\{(x_{3}^{(1)},y_{3}^{(1)},z_{3}^{(1)},w_{3}^{(1)})|x_{3}^{(1)}=y_{3}^{(1)}=0\}{\rm : \rm}
\end{equation*}
$$
{x_{3}^{(2)}}=\frac{x_{3}^{(1)}}{y_{3}^{(1)}} \;, \;\;\; {y_{3}^{(2)}}=y_{3}^{(1)} \;, \;\;\; {z_{3}^{(2)}}=z_{3}^{(1)} \;, \;\;\; {w_{3}^{(2)}}=w_{3}^{(1)}.
$$

{\bf Step 3:} We blow up along the surface
\begin{equation*}
\{(x_{3}^{(2)},y_{3}^{(2)},z_{3}^{(2)},w_{3}^{(2)})|x_{3}^{(2)}-\alpha_3=y_{3}^{(2)}=0\}{\rm : \rm}
\end{equation*}
$$
{x_{3}^{(3)}}=\frac{x_{3}^{(2)}-\alpha_3}{y_{3}^{(2)}} \;, \;\;\; {y_{3}^{(3)}}=y_{3}^{(2)} \;, \;\;\; {z_{3}^{(3)}}=z_{3}^{(2)} \;, \;\;\; {w_{3}^{(3)}}=w_{3}^{(2)}.
$$
We have resolved the accessible singular point $P_3$.

By choosing a new coordinate system as
$$
(x_3,y_3,z_3,w_3)=(-x_{3}^{(3)},y_{3}^{(3)},z_{3}^{(3)},w_{3}^{(3)}),
$$
we can obtain the coordinate system $(x_3,y_3,z_3,w_3)$ given in Theorem \ref{th;holoG}.

Next, we give an explicit resolution process for the accessible singular point $P_6$ by giving a convenient coordinate system at each step.

By the following steps, we can resolve the accessible singular point $P_6$.

{\bf Step 0:} We take the coordinate system centered at $P_6${\rm : \rm}
$$
{x_{6}^{(0)}}=X_3 \;, \;\;\; {y_{6}^{(0)}}=Y_3 \;, \;\;\; {z_{6}^{(0)}}=Z_3 \;, \;\;\; {w_{6}^{(0)}}=W_3-t/s.
$$

{\bf Step 1:} We blow up at the point $P_6${\rm : \rm}
$$
{x_{6}^{(1)}}=\frac{x_{6}^{(0)}}{y_{6}^{(0)}} \;, \;\;\; {y_{6}^{(1)}}=y_{6}^{(0)} \;, \;\;\; {z_{6}^{(1)}}=\frac{z_{6}^{(0)}}{y_{6}^{(0)}} \;, \;\;\; {w_{6}^{(1)}}=\frac{w_{6}^{(0)}}{y_{6}^{(0)}}.
$$

{\bf Step 2:} We blow up along the surface
\begin{equation*}
\{(x_{6}^{(1)},y_{6}^{(1)},z_{6}^{(1)},w_{6}^{(1)})|x_{6}^{(1)}+tz_{6}^{(1)}/s-t=y_{6}^{(1)}=0\}{\rm : \rm}
\end{equation*}
$$
{x_{6}^{(2)}}=\frac{x_{6}^{(1)}+tz_{6}^{(1)}/s-t}{y_{6}^{(1)}} \;, \;\;\; {y_{6}^{(2)}}=y_{6}^{(1)} \;, \;\;\; {z_{6}^{(2)}}=z_{6}^{(1)} \;, \;\;\; {w_{6}^{(2)}}=w_{6}^{(1)}.
$$

{\bf Step 3:} We blow up along the surface
\begin{equation*}
\{(x_{6}^{(2)},y_{6}^{(2)},z_{6}^{(2)},w_{6}^{(2)})|x_{6}^{(2)}-\alpha_6=y_{6}^{(2)}=0\}{\rm : \rm}
\end{equation*}
$$
{x_{6}^{(3)}}=\frac{x_{6}^{(2)}-\alpha_6}{y_{6}^{(2)}} \;, \;\;\; {y_{6}^{(3)}}=y_{6}^{(2)} \;, \;\;\; {z_{6}^{(3)}}=z_{6}^{(2)} \;, \;\;\; {w_{6}^{(3)}}=w_{6}^{(2)}.
$$
We have resolved the accessible singular point $P_6$.

By choosing a new coordinate system as
$$
(x_6,y_6,z_6,w_6)=(-x_{6}^{(3)},y_{6}^{(3)},z_{6}^{(3)},w_{6}^{(3)}),
$$
we can obtain the coordinate system $(x_6,y_6,z_6,w_6)$ given in Theorem \ref{th;holoG}.

For the cases of $P_4,P_5$, the proof is similar.

The proof has thus been completed.  \qed

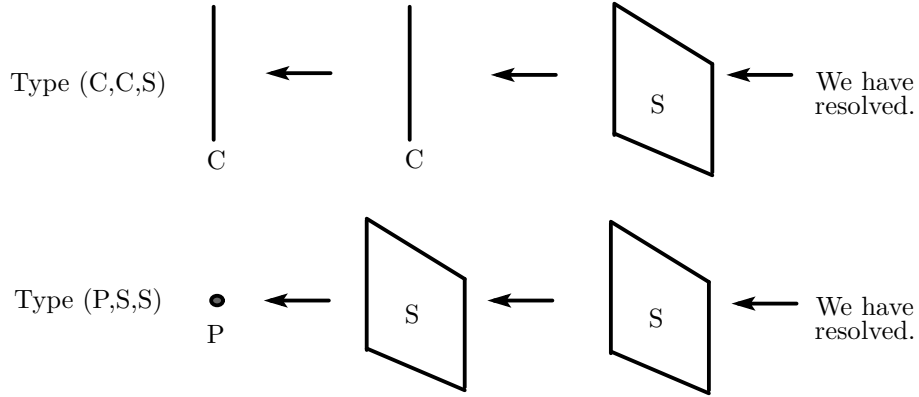
\begin{figure}
\unitlength 0.1in
\begin{picture}(42.15,20.41)(3.40,-23.51)
\put(3.4000,-7.6900){\makebox(0,0)[lb]{Type (C,C,S)}}%
%
\special{pn 20}%
\special{pa 1405 326}%
\special{pa 1405 1023}%
\special{fp}%
%
\special{pn 20}%
\special{pa 2017 674}%
\special{pa 1711 674}%
\special{fp}%
\special{sh 1}%
\special{pa 1711 674}%
\special{pa 1778 694}%
\special{pa 1764 674}%
\special{pa 1778 654}%
\special{pa 1711 674}%
\special{fp}%
%
\special{pn 20}%
\special{pa 2431 326}%
\special{pa 2431 1023}%
\special{fp}%
%
\special{pn 20}%
\special{pa 3196 683}%
\special{pa 2890 683}%
\special{fp}%
\special{sh 1}%
\special{pa 2890 683}%
\special{pa 2957 703}%
\special{pa 2943 683}%
\special{pa 2957 663}%
\special{pa 2890 683}%
\special{fp}%
%
\special{pn 20}%
\special{pa 3502 310}%
\special{pa 4015 626}%
\special{fp}%
\special{pa 3502 318}%
\special{pa 3502 998}%
\special{fp}%
\special{pa 4015 626}%
\special{pa 4015 1193}%
\special{fp}%
\special{pa 3502 990}%
\special{pa 4015 1209}%
\special{fp}%
%
\special{pn 20}%
\special{pa 4429 683}%
\special{pa 4123 683}%
\special{fp}%
\special{sh 1}%
\special{pa 4123 683}%
\special{pa 4190 703}%
\special{pa 4176 683}%
\special{pa 4190 663}%
\special{pa 4123 683}%
\special{fp}%
\put(45.5500,-7.6400){\makebox(0,0)[lb]{We have}}%
\put(3.6000,-19.0300){\makebox(0,0)[lb]{Type (P,S,S)}}%
%
\special{pn 20}%
\special{sh 0.600}%
\special{ar 1423 1865 31 25  0.0000000 6.2831853}%
%
\special{pn 20}%
\special{pa 2008 1865}%
\special{pa 1702 1865}%
\special{fp}%
\special{sh 1}%
\special{pa 1702 1865}%
\special{pa 1769 1885}%
\special{pa 1755 1865}%
\special{pa 1769 1845}%
\special{pa 1702 1865}%
\special{fp}%
%
\special{pn 20}%
\special{pa 2206 1436}%
\special{pa 2719 1752}%
\special{fp}%
\special{pa 2206 1444}%
\special{pa 2206 2124}%
\special{fp}%
\special{pa 2719 1752}%
\special{pa 2719 2319}%
\special{fp}%
\special{pa 2206 2116}%
\special{pa 2719 2335}%
\special{fp}%
%
\special{pn 20}%
\special{pa 3178 1865}%
\special{pa 2872 1865}%
\special{fp}%
\special{sh 1}%
\special{pa 2872 1865}%
\special{pa 2939 1885}%
\special{pa 2925 1865}%
\special{pa 2939 1845}%
\special{pa 2872 1865}%
\special{fp}%
%
\special{pn 20}%
\special{pa 3484 1452}%
\special{pa 3997 1768}%
\special{fp}%
\special{pa 3484 1460}%
\special{pa 3484 2141}%
\special{fp}%
\special{pa 3997 1768}%
\special{pa 3997 2335}%
\special{fp}%
\special{pa 3484 2132}%
\special{pa 3997 2351}%
\special{fp}%
%
\special{pn 20}%
\special{pa 4456 1881}%
\special{pa 4150 1881}%
\special{fp}%
\special{sh 1}%
\special{pa 4150 1881}%
\special{pa 4217 1901}%
\special{pa 4203 1881}%
\special{pa 4217 1861}%
\special{pa 4150 1881}%
\special{fp}%
\put(45.5500,-19.3800){\makebox(0,0)[lb]{We have}}%
\put(13.6900,-11.7700){\makebox(0,0)[lb]{C}}%
\put(24.0400,-11.6900){\makebox(0,0)[lb]{C}}%
\put(36.9100,-8.8500){\makebox(0,0)[lb]{S}}%
\put(13.6900,-20.8400){\makebox(0,0)[lb]{P}}%
\put(24.0400,-19.7900){\makebox(0,0)[lb]{S}}%
\put(36.8200,-20.0300){\makebox(0,0)[lb]{S}}%
\put(45.5000,-8.8600){\makebox(0,0)[lb]{resolved.}}%
\put(45.5000,-20.6500){\makebox(0,0)[lb]{resolved.}}%
\end{picture}%
\label{fig:G4}
\caption{This figure denotes resolution process of the accessible singularities. We denote that P is a point, C is a curve and S is a surface.}
\end{figure}

\section{New expression of the system \eqref{1}}
In the case of holomorphy conditions $(x_i,y_i,z_i,w_i) \ (i=1,2,\dots,6)$, each coordinate system is classified as follows:

\begin{center}
\begin{tabular}{|c|c|} \hline 
Type (C,C,S) & $(x_i,y_i,z_i,w_i) \ (i=1,2)$    \\ \hline 
Type (P,S,S) & $(x_i,y_i,z_i,w_i) \ (i=3,4,5,6)$    \\ \hline 
\end{tabular}
\end{center}
We denote that P is a point, C is a curve and S is a surface (see Figure 2). These properties suggest the possibility that there suggests a procedure for searching for other versions with different types of accessible singularities from the system \eqref{1}. Now, we consider the following problem.

\begin{Problem}
Can we find a polynomial Hamiltonian system with Hamiltonians $H_i \in {\Bbb C}(t,s)[q_1,p_1,q_2,p_2] \ (i=1,2)$ with different types of accessible singularities from the system \eqref{1}?
\end{Problem}

To answer this, we construct the polynomial Hamiltonian system \eqref{uraS} with different type of holomorphy conditions of the system \eqref{1}.

\begin{thm}\label{holoURA}
Let us consider a polynomial Hamiltonian system with Hamiltonians $H_i \in {\Bbb C}(t,s)[q_1,p_1,q_2,p_2] \ (i=1,2)$. We assume that

$(A1)$ $deg(H_i)=5$ with respect to $q_1,p_1,q_2,p_2$.

$(A2)$ This system becomes again a polynomial Hamiltonian system in each coordinate $r'_j, \ j=1,2,\dots,6${\rm : \rm}
\begin{align}
\begin{split}
r'_1:&x_1=\frac{1}{q_1}, \ y_1=-q_1(q_1p_1+\alpha_2), \ z_1=q_2, \ w_1=p_2,\\
r'_2:&x_2=\frac{1}{q_1}, \ y_2=-q_1(q_1p_1+q_2p_2-\alpha_1), \ z_2=\frac{q_2}{q_1}, \ w_2=q_1p_2,\\
r'_3:&x_3=q_1,\ y_3=p_1,\ z_3=\frac{1}{q_2},\ w_3=-(q_2p_2+\alpha_4)q_2,\\
r'_4:&x_4=-(q_1p_1+q_2p_2-(\alpha_1+\alpha_3))p_1, \ y_4=\frac{1}{p_1}, \ z_4=q_2p_1, \ w_4=\frac{p_2}{p_1},\\
r'_5:&x_5=-((q_1-1)p_1+(q_2-1)p_2-(\alpha_1+\alpha_5))p_1, \ y_5=\frac{1}{p_1},\\
&z_5=(q_2-1)p_1, \ w_5=\frac{p_2}{p_1},\\
r'_6:&x_6=-((q_1-t)p_1+(q_2-s)p_2-(\alpha_1+\alpha_6))p_1, \ y_6=\frac{1}{p_1},\\
&z_6=(q_2-s)p_1, \ w_6=\frac{p_2}{p_1}.
\end{split}
\end{align}
Then such a system coincides with the system \eqref{uraS}.
\end{thm}
We remark that each transformation of each coordinate $r'_i \ (i=1,2,\dots,6)$ is birational and symplectic.

\begin{prop} 
On each coordinate $r'_j, \ j=1,2,\dots,6$, the  Hamiltonians $H'_{j1}$ and $H'_{j2}$ on $U_j \times B$ are expressed as a polynomial in $x_j,y_j,z_j,w_j$ and a rational function in $t$ and $s$, and satisfy the following conditions{\rm: \rm}
\begin{align}
\begin{split}
&dq_1 \wedge dp_1 +dq_2 \wedge dp_2 - dH'_1 \wedge dt- dH'_2 \wedge ds\\
&=dx_j \wedge dy_j +dz_j \wedge dw_j - dH'_{j1} \wedge dt- dH'_{j2} \wedge ds \quad (j=1,2,\dots,5),\\
&dq_1 \wedge dp_1 +dq_2 \wedge dp_2 - d(H'_1-p_1) \wedge dt- d(H'_2-p_2) \wedge ds\\
&=dx_6 \wedge dy_6 +dz_6 \wedge dw_6 - dH'_{61} \wedge dt- dH'_{62} \wedge ds.
\end{split}
\end{align}
\end{prop}

\noindent
{\bf Proof of Theorem \ref{holoURA}.}
For the transformations $r'_i \ (i=1,2, \dots,5)$ we can repeat the same way given in Theorem \ref{th;holoG}.

We must be careful of the case of the transformation $r'_6$. The relation between the coordinate system $(x_6,y_6,z_6,w_6)$ and the coordinate system $(q_1,p_1,q_2,p_2)$ is given by
\begin{align}
\begin{split}
&dx_6 \wedge dy_6 +dz_6 \wedge dw_6\\
&=dq_1 \wedge dp_1 +dq_2 \wedge dp_2+dp_1 \wedge dt+dp_2 \wedge ds.
\end{split}
\end{align}
Resolving the coordinate $r_6$ in the variables $q_1,p_1,q_2,p_2$, we obtain
\begin{equation}\tag*{($R_6$)}
q_1=t-x_6y_6^2-y_6z_6w_6+(\alpha_1+\alpha_6)y_6, \quad p_1=1/y_6, \quad q_2=y_6z_6+s, \quad p_2=w_6/y_6.
\end{equation}
In this case, we must consider the polynomiality for $R_6(H_1-p_1)$ and $R_6(H_2-p_2)$, respectively.

In this way, we can obtain the Hamiltonians $H'_1,H'_2$. \qed

In the case of holomorphy conditions $(x_i,y_i,z_i,w_i) \ (i=1,2,\dots,6)$, each coordinate system is classified as follows:

\begin{center}
\begin{tabular}{|c|c|} \hline 
Type (C,C,S) & $(x_i,y_i,z_i,w_i) \ (i=2,4,5,6)$    \\ \hline 
Type (P,S,S) & $(x_i,y_i,z_i,w_i) \ (i=1,3)$    \\ \hline 
\end{tabular}
\end{center}

Moreover, we show that the system \eqref{1} is equivalent to the system \eqref{uraS} by giving an explicit rational and symplectic transformation. This transformation can be considered as a generalization of {\it Okamoto-transformation \it}of the sixth Painlev\'e system.

\begin{thm}\label{th:oka}
By using the rational and symplectic transformation \eqref{extensionoka}, the system \eqref{1} coincides with \eqref{uraS}.
\end{thm}

{\bf Proof of Theorem \ref{th:oka}}
Set
\begin{equation*}
X:=q_1+\frac{q_2p_2+\alpha_1}{p_1}, \quad Y:=p_1, \quad Z:=\frac{p_2}{p_1}, \quad W:=-q_2p_1 \quad T:=t, \quad S:=\frac{t}{s}.
\end{equation*}
By resolving in $q_1,p_1,q_2,p_2,t,s$, we obtain
\begin{equation}\tag*{$\tilde{S}:$}
q_1=X+\frac{ZW-\alpha_1}{Y}, \quad p_1=Y, \quad q_2=-\frac{W}{Y}, \quad p_2=YZ, \quad t=T, \quad s=\frac{T}{S}.
\end{equation}
Applying the transformations in $t,s$ and the transformation of the symplectic 2-form:
\begin{align}
\begin{split}
&dX \wedge dY +dZ \wedge dW=dq_1 \wedge dp_1 +dq_2 \wedge dp_2,\\
&dt=dT,\\
&ds=\frac{1}{S}dT-\frac{T}{S^2}dS,
\end{split}
\end{align}
we obtain the polynomial Hamiltonians $\tilde{S}(H_1+\frac{1}{S}H_2),\tilde{S}(-\frac{T}{S^2}H_2)$, which satisfy the following relations:
$$
\tilde{S}(H_1+\frac{1}{S}H_2)=H'_1, \quad \tilde{S}(-\frac{T}{S^2}H_2)=H'_2.
$$
This completes the proof. \qed

We also consider the inverse transformation of \eqref{extensionoka}.
\begin{Lemma}\label{th:inoka}
The transformation
\begin{equation}
R:(q_1,p_1,q_2,p_2,t,s) \rightarrow \left(q_1+\frac{q_2p_2-\alpha_1}{p_1},p_1,-\frac{p_2}{p_1},q_2p_1,t,\frac{t}{s}\right)
\end{equation}
satisfy the following relations\rm{:\rm}
\begin{equation}
R \circ S=1, \quad S \circ R=1.
\end{equation}
\end{Lemma}

Composing such transformations, we can make a B{\"a}cklund transformation.
\begin{prop}
By using the transformations $S$ and
\begin{equation}
R':(q_1,p_1,q_2,p_2,t,s) \rightarrow \left(-\frac{p_1}{p_2},q_1p_2,q_2+\frac{q_1p_1-\alpha_1}{p_2},p_2,\frac{s}{t},s\right),
\end{equation}
we make the B{\"a}cklund transformation $R' \circ S$ of \eqref{1}\rm{:\rm}
\begin{align*}
R' \circ S:&(q_1,p_1,q_2,p_2,t,s;\alpha_1,\dots,\alpha_6) \rightarrow\\
&\left(\frac{1}{q_2},-(q_1p_1+q_2p_2+\alpha_1)q_2,-\frac{q_1}{q_2},-p_1q_2,\frac{1}{s},\frac{t}{s};\alpha_1,\alpha_4,\alpha_2,\alpha_3,\alpha_5,\alpha_6\right).
\end{align*}
\end{prop}

\begin{prop}
By using the transformations $R$ and
\begin{align*}
N:(q_1,p_1,q_2,p_2t,s) \rightarrow (&\frac{p_1}{q_1p_1+q_2p_2+\alpha_1},-q_1(q_1p_1+q_2p_2+\alpha_1),\\
&\frac{p_2}{q_1p_1+q_2p_2+\alpha_1},-q_2(q_1p_1+q_2p_2+\alpha_1),\frac{1}{t},\frac{1}{s}),
\end{align*}
we make the B{\"a}cklund transformation $R \circ N$ of \eqref{uraS}\rm{:\rm}
\begin{align*}
R \circ N:&(q_1,p_1,q_2,p_2,t,s;\alpha_1.\dots,\alpha_6) \rightarrow\\
&\left(\frac{1}{q_1},-(q_1p_1+q_2p_2+\alpha_1)q_1,-\frac{q_2}{q_1},-p_2q_1,\frac{1}{t},\frac{s}{t};\alpha_1,\alpha_3,\alpha_2,\alpha_4,\alpha_5,\alpha_6\right).
\end{align*}
\end{prop}

Now, we consider the following problem.
\begin{Problem}
For the system \eqref{uraS}, can we find such a $S$-transformation differ-\\
ent from the transformation $R$?
\end{Problem}

\begin{figure}
\unitlength 0.1in
\begin{picture}(36.79,27.48)(19.60,-28.45)
\put(33.1000,-22.3000){\makebox(0,0)[lb]{degree 5}}%
\put(31.4000,-25.4000){\makebox(0,0)[lb]{degree 6}}%
%
\special{pn 20}%
\special{pa 3300 2140}%
\special{pa 3268 2126}%
\special{pa 3236 2112}%
\special{pa 3204 2097}%
\special{pa 3173 2083}%
\special{pa 3142 2068}%
\special{pa 3112 2052}%
\special{pa 3083 2036}%
\special{pa 3055 2019}%
\special{pa 3027 2002}%
\special{pa 3002 1983}%
\special{pa 2977 1964}%
\special{pa 2955 1943}%
\special{pa 2934 1922}%
\special{pa 2914 1899}%
\special{pa 2897 1875}%
\special{pa 2882 1849}%
\special{pa 2869 1822}%
\special{pa 2859 1794}%
\special{pa 2850 1764}%
\special{pa 2843 1734}%
\special{pa 2839 1702}%
\special{pa 2836 1670}%
\special{pa 2835 1637}%
\special{pa 2836 1604}%
\special{pa 2839 1570}%
\special{pa 2843 1536}%
\special{pa 2850 1502}%
\special{pa 2857 1468}%
\special{pa 2867 1434}%
\special{pa 2878 1400}%
\special{pa 2890 1366}%
\special{pa 2904 1334}%
\special{pa 2919 1301}%
\special{pa 2936 1270}%
\special{pa 2953 1240}%
\special{pa 2972 1210}%
\special{pa 2993 1182}%
\special{pa 3014 1155}%
\special{pa 3036 1129}%
\special{pa 3060 1105}%
\special{pa 3084 1081}%
\special{pa 3110 1059}%
\special{pa 3136 1038}%
\special{pa 3163 1019}%
\special{pa 3191 1000}%
\special{pa 3220 983}%
\special{pa 3249 967}%
\special{pa 3279 952}%
\special{pa 3310 939}%
\special{pa 3341 926}%
\special{pa 3373 915}%
\special{pa 3406 905}%
\special{pa 3439 895}%
\special{pa 3472 887}%
\special{pa 3506 881}%
\special{pa 3540 875}%
\special{pa 3575 870}%
\special{pa 3609 867}%
\special{pa 3644 864}%
\special{pa 3679 863}%
\special{pa 3714 862}%
\special{pa 3750 863}%
\special{pa 3785 865}%
\special{pa 3820 868}%
\special{pa 3856 871}%
\special{pa 3891 876}%
\special{pa 3926 882}%
\special{pa 3961 889}%
\special{pa 3995 896}%
\special{pa 4030 905}%
\special{pa 4064 915}%
\special{pa 4097 926}%
\special{pa 4131 937}%
\special{pa 4164 950}%
\special{pa 4196 963}%
\special{pa 4228 978}%
\special{pa 4259 993}%
\special{pa 4290 1009}%
\special{pa 4320 1027}%
\special{pa 4349 1045}%
\special{pa 4378 1064}%
\special{pa 4405 1083}%
\special{pa 4432 1104}%
\special{pa 4457 1126}%
\special{pa 4482 1148}%
\special{pa 4505 1172}%
\special{pa 4527 1196}%
\special{pa 4548 1221}%
\special{pa 4567 1247}%
\special{pa 4585 1274}%
\special{pa 4602 1301}%
\special{pa 4617 1330}%
\special{pa 4630 1359}%
\special{pa 4641 1389}%
\special{pa 4651 1420}%
\special{pa 4659 1452}%
\special{pa 4665 1484}%
\special{pa 4669 1518}%
\special{pa 4671 1551}%
\special{pa 4672 1585}%
\special{pa 4670 1619}%
\special{pa 4667 1653}%
\special{pa 4662 1687}%
\special{pa 4655 1720}%
\special{pa 4646 1753}%
\special{pa 4636 1785}%
\special{pa 4624 1816}%
\special{pa 4610 1846}%
\special{pa 4594 1874}%
\special{pa 4577 1901}%
\special{pa 4558 1926}%
\special{pa 4538 1949}%
\special{pa 4516 1970}%
\special{pa 4492 1989}%
\special{pa 4467 2007}%
\special{pa 4441 2022}%
\special{pa 4413 2036}%
\special{pa 4385 2049}%
\special{pa 4355 2060}%
\special{pa 4325 2070}%
\special{pa 4293 2078}%
\special{pa 4261 2086}%
\special{pa 4228 2093}%
\special{pa 4195 2099}%
\special{pa 4161 2105}%
\special{pa 4126 2110}%
\special{pa 4092 2115}%
\special{pa 4057 2119}%
\special{pa 4029 2123}%
\special{sp}%
%
\special{pn 20}%
\special{pa 3030 2460}%
\special{pa 2999 2444}%
\special{pa 2969 2428}%
\special{pa 2938 2411}%
\special{pa 2908 2395}%
\special{pa 2878 2378}%
\special{pa 2848 2362}%
\special{pa 2819 2345}%
\special{pa 2790 2327}%
\special{pa 2761 2310}%
\special{pa 2734 2292}%
\special{pa 2706 2273}%
\special{pa 2680 2255}%
\special{pa 2654 2235}%
\special{pa 2629 2216}%
\special{pa 2605 2195}%
\special{pa 2581 2174}%
\special{pa 2559 2153}%
\special{pa 2538 2131}%
\special{pa 2518 2108}%
\special{pa 2499 2084}%
\special{pa 2481 2060}%
\special{pa 2464 2035}%
\special{pa 2449 2008}%
\special{pa 2435 1981}%
\special{pa 2423 1954}%
\special{pa 2412 1925}%
\special{pa 2402 1896}%
\special{pa 2394 1866}%
\special{pa 2387 1835}%
\special{pa 2381 1804}%
\special{pa 2377 1772}%
\special{pa 2374 1739}%
\special{pa 2371 1707}%
\special{pa 2371 1674}%
\special{pa 2371 1640}%
\special{pa 2372 1607}%
\special{pa 2375 1573}%
\special{pa 2379 1539}%
\special{pa 2383 1505}%
\special{pa 2389 1470}%
\special{pa 2396 1436}%
\special{pa 2404 1402}%
\special{pa 2412 1368}%
\special{pa 2422 1334}%
\special{pa 2433 1300}%
\special{pa 2445 1267}%
\special{pa 2457 1234}%
\special{pa 2471 1201}%
\special{pa 2485 1168}%
\special{pa 2500 1136}%
\special{pa 2516 1105}%
\special{pa 2533 1074}%
\special{pa 2550 1044}%
\special{pa 2568 1014}%
\special{pa 2587 986}%
\special{pa 2607 958}%
\special{pa 2627 930}%
\special{pa 2648 904}%
\special{pa 2670 879}%
\special{pa 2692 854}%
\special{pa 2715 831}%
\special{pa 2738 808}%
\special{pa 2762 787}%
\special{pa 2787 766}%
\special{pa 2812 746}%
\special{pa 2838 727}%
\special{pa 2864 709}%
\special{pa 2890 691}%
\special{pa 2918 674}%
\special{pa 2945 659}%
\special{pa 2973 643}%
\special{pa 3002 629}%
\special{pa 3031 615}%
\special{pa 3060 602}%
\special{pa 3090 590}%
\special{pa 3120 578}%
\special{pa 3151 567}%
\special{pa 3181 557}%
\special{pa 3213 547}%
\special{pa 3244 538}%
\special{pa 3276 529}%
\special{pa 3308 521}%
\special{pa 3341 513}%
\special{pa 3374 506}%
\special{pa 3407 500}%
\special{pa 3440 493}%
\special{pa 3473 488}%
\special{pa 3507 482}%
\special{pa 3541 478}%
\special{pa 3575 473}%
\special{pa 3609 469}%
\special{pa 3644 465}%
\special{pa 3679 462}%
\special{pa 3713 459}%
\special{pa 3748 456}%
\special{pa 3783 454}%
\special{pa 3818 452}%
\special{pa 3853 451}%
\special{pa 3888 450}%
\special{pa 3923 449}%
\special{pa 3958 449}%
\special{pa 3993 450}%
\special{pa 4028 451}%
\special{pa 4063 452}%
\special{pa 4097 454}%
\special{pa 4132 456}%
\special{pa 4166 459}%
\special{pa 4201 463}%
\special{pa 4235 466}%
\special{pa 4269 471}%
\special{pa 4302 476}%
\special{pa 4336 482}%
\special{pa 4369 488}%
\special{pa 4402 495}%
\special{pa 4435 502}%
\special{pa 4467 510}%
\special{pa 4499 519}%
\special{pa 4530 528}%
\special{pa 4562 538}%
\special{pa 4593 548}%
\special{pa 4623 559}%
\special{pa 4653 571}%
\special{pa 4682 584}%
\special{pa 4711 597}%
\special{pa 4740 611}%
\special{pa 4768 625}%
\special{pa 4795 641}%
\special{pa 4822 657}%
\special{pa 4849 673}%
\special{pa 4874 691}%
\special{pa 4900 709}%
\special{pa 4924 728}%
\special{pa 4948 748}%
\special{pa 4971 769}%
\special{pa 4993 790}%
\special{pa 5015 813}%
\special{pa 5036 836}%
\special{pa 5056 860}%
\special{pa 5076 884}%
\special{pa 5095 910}%
\special{pa 5112 936}%
\special{pa 5130 963}%
\special{pa 5146 990}%
\special{pa 5161 1019}%
\special{pa 5176 1047}%
\special{pa 5190 1077}%
\special{pa 5203 1107}%
\special{pa 5215 1137}%
\special{pa 5227 1168}%
\special{pa 5237 1199}%
\special{pa 5247 1231}%
\special{pa 5256 1263}%
\special{pa 5263 1296}%
\special{pa 5271 1328}%
\special{pa 5277 1361}%
\special{pa 5282 1395}%
\special{pa 5286 1428}%
\special{pa 5290 1462}%
\special{pa 5292 1495}%
\special{pa 5294 1529}%
\special{pa 5294 1563}%
\special{pa 5294 1597}%
\special{pa 5293 1631}%
\special{pa 5291 1665}%
\special{pa 5288 1698}%
\special{pa 5283 1732}%
\special{pa 5278 1765}%
\special{pa 5272 1799}%
\special{pa 5265 1832}%
\special{pa 5257 1865}%
\special{pa 5248 1897}%
\special{pa 5238 1929}%
\special{pa 5227 1961}%
\special{pa 5215 1992}%
\special{pa 5202 2022}%
\special{pa 5188 2052}%
\special{pa 5173 2081}%
\special{pa 5158 2109}%
\special{pa 5141 2136}%
\special{pa 5123 2163}%
\special{pa 5105 2188}%
\special{pa 5085 2212}%
\special{pa 5065 2235}%
\special{pa 5044 2256}%
\special{pa 5022 2276}%
\special{pa 4999 2295}%
\special{pa 4975 2313}%
\special{pa 4950 2330}%
\special{pa 4925 2345}%
\special{pa 4899 2359}%
\special{pa 4872 2372}%
\special{pa 4844 2385}%
\special{pa 4816 2396}%
\special{pa 4787 2406}%
\special{pa 4758 2415}%
\special{pa 4728 2424}%
\special{pa 4697 2431}%
\special{pa 4666 2438}%
\special{pa 4634 2444}%
\special{pa 4602 2450}%
\special{fp}%
\special{pa 4602 2450}%
\special{pa 4569 2454}%
\special{pa 4536 2458}%
\special{pa 4503 2462}%
\special{pa 4469 2465}%
\special{pa 4435 2467}%
\special{pa 4400 2469}%
\special{pa 4366 2470}%
\special{pa 4331 2472}%
\special{pa 4295 2472}%
\special{pa 4260 2473}%
\special{pa 4224 2473}%
\special{pa 4188 2473}%
\special{pa 4152 2472}%
\special{pa 4116 2472}%
\special{pa 4080 2471}%
\special{pa 4044 2471}%
\special{pa 4010 2470}%
\special{sp}%
\put(29.2000,-16.2000){\makebox(0,0)[lb]{$(H_1,H_2)$}}%
\put(38.7000,-16.2000){\makebox(0,0)[lb]{$(H'_1,H'_2)$}}%
%
\special{pn 20}%
\special{pa 3530 1560}%
\special{pa 3860 1560}%
\special{fp}%
\special{sh 1}%
\special{pa 3860 1560}%
\special{pa 3793 1540}%
\special{pa 3807 1560}%
\special{pa 3793 1580}%
\special{pa 3860 1560}%
\special{fp}%
%
\special{pn 20}%
\special{pa 4190 1670}%
\special{pa 4700 1980}%
\special{fp}%
\special{sh 1}%
\special{pa 4700 1980}%
\special{pa 4653 1928}%
\special{pa 4654 1952}%
\special{pa 4633 1962}%
\special{pa 4700 1980}%
\special{fp}%
\put(44.4000,-21.7000){\makebox(0,0)[lb]{$(H''_1,H''_2)$}}%
\put(36.0000,-14.7000){\makebox(0,0)[lb]{$S$}}%
\put(44.1000,-18.1000){\makebox(0,0)[lb]{$S_1$}}%
%
\special{pn 20}%
\special{pa 2980 2840}%
\special{pa 2949 2825}%
\special{pa 2918 2809}%
\special{pa 2887 2794}%
\special{pa 2856 2779}%
\special{pa 2825 2763}%
\special{pa 2794 2748}%
\special{pa 2763 2732}%
\special{pa 2733 2716}%
\special{pa 2703 2700}%
\special{pa 2673 2684}%
\special{pa 2643 2668}%
\special{pa 2614 2651}%
\special{pa 2585 2634}%
\special{pa 2556 2617}%
\special{pa 2528 2600}%
\special{pa 2500 2582}%
\special{pa 2472 2564}%
\special{pa 2445 2546}%
\special{pa 2419 2528}%
\special{pa 2392 2509}%
\special{pa 2367 2489}%
\special{pa 2342 2470}%
\special{pa 2318 2450}%
\special{pa 2294 2429}%
\special{pa 2271 2408}%
\special{pa 2248 2387}%
\special{pa 2226 2365}%
\special{pa 2205 2342}%
\special{pa 2185 2319}%
\special{pa 2165 2296}%
\special{pa 2147 2271}%
\special{pa 2129 2247}%
\special{pa 2112 2221}%
\special{pa 2095 2195}%
\special{pa 2080 2169}%
\special{pa 2066 2141}%
\special{pa 2052 2113}%
\special{pa 2040 2085}%
\special{pa 2028 2056}%
\special{pa 2018 2026}%
\special{pa 2008 1996}%
\special{pa 1999 1965}%
\special{pa 1991 1934}%
\special{pa 1984 1902}%
\special{pa 1978 1870}%
\special{pa 1973 1838}%
\special{pa 1969 1805}%
\special{pa 1965 1772}%
\special{pa 1963 1739}%
\special{pa 1961 1705}%
\special{pa 1960 1671}%
\special{pa 1960 1638}%
\special{pa 1961 1604}%
\special{pa 1962 1570}%
\special{pa 1965 1536}%
\special{pa 1968 1502}%
\special{pa 1972 1467}%
\special{pa 1977 1433}%
\special{pa 1983 1400}%
\special{pa 1990 1366}%
\special{pa 1997 1332}%
\special{pa 2005 1299}%
\special{pa 2014 1266}%
\special{pa 2023 1233}%
\special{pa 2034 1200}%
\special{pa 2045 1168}%
\special{pa 2057 1136}%
\special{pa 2069 1104}%
\special{pa 2083 1073}%
\special{pa 2097 1043}%
\special{pa 2112 1013}%
\special{pa 2127 983}%
\special{pa 2143 954}%
\special{pa 2160 926}%
\special{pa 2178 898}%
\special{pa 2196 871}%
\special{pa 2215 845}%
\special{pa 2235 819}%
\special{pa 2255 794}%
\special{pa 2276 769}%
\special{pa 2297 745}%
\special{pa 2319 722}%
\special{pa 2342 699}%
\special{pa 2365 677}%
\special{pa 2388 655}%
\special{pa 2412 634}%
\special{pa 2437 614}%
\special{pa 2462 594}%
\special{pa 2488 574}%
\special{pa 2514 555}%
\special{pa 2540 537}%
\special{pa 2567 519}%
\special{pa 2594 501}%
\special{pa 2622 484}%
\special{pa 2650 467}%
\special{pa 2678 451}%
\special{pa 2707 435}%
\special{pa 2736 420}%
\special{pa 2765 405}%
\special{pa 2795 390}%
\special{pa 2825 376}%
\special{pa 2855 362}%
\special{pa 2885 349}%
\special{pa 2916 336}%
\special{pa 2946 323}%
\special{pa 2977 310}%
\special{pa 3008 298}%
\special{pa 3040 286}%
\special{pa 3071 275}%
\special{pa 3102 264}%
\special{pa 3134 253}%
\special{pa 3166 242}%
\special{pa 3197 232}%
\special{pa 3229 222}%
\special{pa 3261 212}%
\special{pa 3293 203}%
\special{pa 3325 194}%
\special{pa 3358 185}%
\special{pa 3390 177}%
\special{pa 3422 169}%
\special{pa 3455 161}%
\special{pa 3487 154}%
\special{pa 3520 148}%
\special{pa 3552 141}%
\special{pa 3585 135}%
\special{pa 3618 130}%
\special{pa 3650 125}%
\special{pa 3683 120}%
\special{pa 3716 116}%
\special{pa 3749 112}%
\special{pa 3781 108}%
\special{pa 3814 105}%
\special{pa 3847 103}%
\special{pa 3880 101}%
\special{pa 3912 99}%
\special{pa 3945 98}%
\special{pa 3978 97}%
\special{pa 4010 97}%
\special{pa 4043 97}%
\special{pa 4075 98}%
\special{pa 4108 100}%
\special{pa 4140 101}%
\special{pa 4173 104}%
\special{pa 4205 107}%
\special{pa 4237 110}%
\special{pa 4269 114}%
\special{pa 4302 118}%
\special{pa 4333 123}%
\special{pa 4365 129}%
\special{pa 4397 135}%
\special{pa 4429 142}%
\special{pa 4460 149}%
\special{pa 4492 157}%
\special{pa 4523 165}%
\special{pa 4554 174}%
\special{pa 4585 184}%
\special{pa 4616 194}%
\special{pa 4647 205}%
\special{pa 4678 216}%
\special{pa 4708 229}%
\special{pa 4738 241}%
\special{pa 4768 255}%
\special{pa 4798 269}%
\special{pa 4828 283}%
\special{pa 4857 299}%
\special{pa 4887 314}%
\special{pa 4916 331}%
\special{pa 4944 348}%
\special{pa 4973 366}%
\special{pa 5001 384}%
\special{pa 5029 403}%
\special{pa 5056 422}%
\special{pa 5083 442}%
\special{pa 5109 463}%
\special{pa 5135 484}%
\special{pa 5161 505}%
\special{pa 5186 527}%
\special{pa 5210 550}%
\special{pa 5234 573}%
\special{pa 5258 596}%
\special{pa 5280 620}%
\special{pa 5302 645}%
\special{pa 5324 669}%
\special{pa 5345 695}%
\special{pa 5365 720}%
\special{pa 5384 746}%
\special{pa 5403 773}%
\special{pa 5421 800}%
\special{pa 5438 827}%
\special{pa 5454 855}%
\special{pa 5469 882}%
\special{pa 5484 911}%
\special{pa 5497 939}%
\special{pa 5510 968}%
\special{pa 5522 998}%
\special{pa 5533 1027}%
\special{pa 5542 1057}%
\special{pa 5551 1087}%
\special{pa 5559 1118}%
\special{pa 5567 1148}%
\special{fp}%
\special{pa 5567 1148}%
\special{pa 5573 1179}%
\special{pa 5579 1210}%
\special{pa 5584 1242}%
\special{pa 5589 1273}%
\special{pa 5593 1305}%
\special{pa 5597 1337}%
\special{pa 5601 1369}%
\special{pa 5604 1401}%
\special{pa 5607 1433}%
\special{pa 5609 1466}%
\special{pa 5612 1498}%
\special{pa 5615 1531}%
\special{pa 5617 1564}%
\special{pa 5620 1597}%
\special{pa 5622 1630}%
\special{pa 5625 1663}%
\special{pa 5628 1696}%
\special{pa 5631 1729}%
\special{pa 5633 1761}%
\special{pa 5635 1794}%
\special{pa 5637 1827}%
\special{pa 5638 1860}%
\special{pa 5639 1892}%
\special{pa 5639 1925}%
\special{pa 5639 1957}%
\special{pa 5638 1989}%
\special{pa 5636 2021}%
\special{pa 5633 2053}%
\special{pa 5629 2084}%
\special{pa 5624 2115}%
\special{pa 5618 2146}%
\special{pa 5610 2177}%
\special{pa 5601 2207}%
\special{pa 5591 2236}%
\special{pa 5580 2266}%
\special{pa 5567 2295}%
\special{pa 5552 2323}%
\special{pa 5536 2351}%
\special{pa 5519 2379}%
\special{pa 5501 2406}%
\special{pa 5482 2432}%
\special{pa 5462 2457}%
\special{pa 5440 2482}%
\special{pa 5418 2507}%
\special{pa 5395 2530}%
\special{pa 5371 2553}%
\special{pa 5346 2575}%
\special{pa 5321 2596}%
\special{pa 5295 2616}%
\special{pa 5268 2636}%
\special{pa 5241 2654}%
\special{pa 5214 2672}%
\special{pa 5186 2688}%
\special{pa 5158 2704}%
\special{pa 5129 2718}%
\special{pa 5101 2732}%
\special{pa 5072 2744}%
\special{pa 5042 2756}%
\special{pa 5013 2767}%
\special{pa 4983 2777}%
\special{pa 4953 2786}%
\special{pa 4922 2794}%
\special{pa 4892 2802}%
\special{pa 4861 2809}%
\special{pa 4830 2815}%
\special{pa 4798 2820}%
\special{pa 4767 2825}%
\special{pa 4735 2830}%
\special{pa 4703 2833}%
\special{pa 4671 2836}%
\special{pa 4639 2839}%
\special{pa 4607 2841}%
\special{pa 4574 2843}%
\special{pa 4541 2844}%
\special{pa 4509 2845}%
\special{pa 4476 2845}%
\special{pa 4443 2845}%
\special{pa 4410 2845}%
\special{pa 4376 2845}%
\special{pa 4343 2844}%
\special{pa 4310 2843}%
\special{pa 4276 2841}%
\special{pa 4243 2840}%
\special{pa 4209 2838}%
\special{pa 4175 2837}%
\special{pa 4142 2835}%
\special{pa 4108 2833}%
\special{pa 4074 2831}%
\special{pa 4060 2830}%
\special{sp}%
\put(31.0000,-29.3000){\makebox(0,0)[lb]{degree 7}}%
\put(45.2000,-26.6000){\makebox(0,0)[lb]{$(H'''_1,H'''_2)$}}%
\put(50.2000,-24.2000){\makebox(0,0)[lb]{$S_2$}}%
%
\special{pn 20}%
\special{pa 4660 2200}%
\special{pa 5020 2470}%
\special{fp}%
\special{sh 1}%
\special{pa 5020 2470}%
\special{pa 4979 2414}%
\special{pa 4977 2438}%
\special{pa 4955 2446}%
\special{pa 5020 2470}%
\special{fp}%
\end{picture}%
\label{fig:Gar1}
\caption{We note that the degree of each $H_i \ (i=1,2)$ is 5 with respect to $q_1,p_1,q_2,p_2$.}
\end{figure}
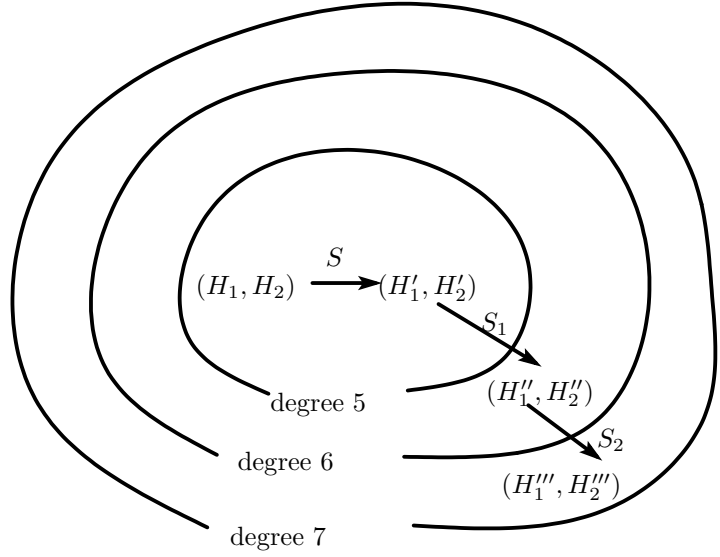

To answer this, we find the following rational and symplectic transformation.
\begin{thm}\label{th:ininoka}
By using the rational and symplectic transformation $S_1$
\begin{equation}
S_1:(q_1,p_1,q_2,p_2,t,s) \rightarrow \left(q_1,p_1+\frac{q_2p_2-\alpha_1-\alpha_3}{q_1},\frac{q_2}{q_1},p_2q_1,t,\frac{s}{t}\right),
\end{equation}
the system \eqref{uraS} is transformed into the Hamiltonian system
\begin{align}\label{urauraS}
\begin{split}
dq_1&=\frac{\partial H_1''}{\partial p_1}dt+\frac{\partial H_2''}{\partial p_1}ds, \quad dp_1=-\frac{\partial H_1''}{\partial q_1}dt-\frac{\partial H_2''}{\partial q_1}ds,\\
dq_2&=\frac{\partial H_1''}{\partial p_2}dt+\frac{\partial H_2''}{\partial p_2}ds, \quad dp_2=-\frac{\partial H_1''}{\partial q_2}dt-\frac{\partial H_2''}{\partial q_2}ds,\\
H_1'' &=H_{VI}(q_1,p_1,t;\alpha_1+\alpha_4+\alpha_6,\alpha_1+\alpha_2,\alpha_3,\alpha_1+\alpha_5,-\alpha_1-\alpha_3)\\
&-\frac{\alpha_4}{t(ts-1)}q_1p_1+\frac{\alpha_3(s-1)}{(t-1)(ts-1)}q_2p_2+\frac{s}{(ts-1)}p_1p_2-\frac{1}{t-1}p_1q_2p_2\\
&+\frac{2(s-1)}{(t-1)(ts-1)}q_1p_1q_2p_2-\frac{q_1q_2(q_1p_1+\alpha_3)\{(ts-1)p_2-(t-1)(q_2p_2+\alpha_4)\}}{t(t-1)(ts-1)},\\
H_2'' &=\pi(H_1''),
\end{split}
\end{align}
where the transformation $\pi$ is explicitly given by
\begin{align}
\begin{split}
&\pi:(q_1,p_1,q_2,p_2,t,s;\alpha_1,\alpha_2, \dots ,\alpha_6)\\
&\rightarrow (q_2,p_2,q_1,p_1,s,t;-\alpha_1-\alpha_2-\alpha_3-\alpha_4,\alpha_2,\alpha_4,\alpha_3,1-\alpha_6,1-\alpha_5).
\end{split}
\end{align}
\end{thm}
We note that each degree of each of Hamiltonians $H''_1,H''_2$ is 6 with respect to $q_1,p_1,q_2,p_2$.

We also remark that on the polynomiality of the transformed Hamiltonian by $S_1$, the change of time-variables $t,s$ of the transformation $S_1$ is not essential.

{\bf Proof of Theorem \ref{th:ininoka}}
Set
\begin{equation*}
X:=q_1, \quad Y:=p_1+\frac{q_2p_2-\alpha_1-\alpha_3}{q_1}, \quad Z:=\frac{q_2}{q_1}, \quad W:=q_1p_2 \quad T:=t, \quad S:=\frac{s}{t}.
\end{equation*}
By resolving in $q_1,p_1,q_2,p_2,t,s$, we obtain
\begin{equation}\tag*{$\tilde{S_1}:$}
q_1=X, \quad p_1=Y-\frac{ZW-\alpha_1-\alpha_3}{X}, \quad q_2=XZ, \quad p_2=\frac{W}{X}, \quad t=T, \quad s=TS.
\end{equation}
Applying the transformations in $t,s$ and the transformation of the symplectic 2-form:
\begin{align}
\begin{split}
&dX \wedge dY +dZ \wedge dW=dq_1 \wedge dp_1 +dq_2 \wedge dp_2,\\
&dt=dT,\\
&ds=SdT+TdS,
\end{split}
\end{align}
we obtain the polynomial Hamiltonians $\tilde{S_1}(H'_1+SH'_2),\tilde{S_1}(TH'_2)$, which satisfy the following relations:
$$
\tilde{S_1}(H'_1+SH'_2)=H''_1, \quad \tilde{S_1}(TH'_2)=H''_2.
$$
This completes the proof. \qed

\begin{thm}\label{holoURAURA}
Let us consider a polynomial Hamiltonian system with Hamiltonians $H_i \in {\Bbb C}(t,s)[q_1,p_1,q_2,p_2] \ (i=1,2)$. We assume that

$(A1)$ $deg(H_i)=6$ with respect to $q_1,p_1,q_2,p_2$.

$(A2)$ This system becomes again a polynomial Hamiltonian system in each coordinate $r''_j, \ j=1,2,\dots,6${\rm : \rm}
\begin{align}
\begin{split}
r''_1:&x_1=-(q_1p_1-q_2p_2+\alpha_1+\alpha_3)p_1, \ y_1=\frac{1}{p_1}, \ z_1=\frac{q_2}{p_1}, \ w_1=p_2p_1,\\
r''_2:&x_2=\frac{1}{q_1}, \ y_2=-q_1(q_1p_1-q_2p_2+\alpha_1+\alpha_2+\alpha_3), \ z_2=q_2q_1, \ w_2=\frac{p_2}{q_1},\\
r''_3:&x_3=\frac{1}{q_1},\ y_3=-(p_1q_1+\alpha_3)q_1,\ z_3=q_2,\ w_3=p_2,\\
r''_4:&x_4=q_1, \ y_4=p_1, \ z_4=\frac{1}{q_2}, \ w_4=-(p_2q_2+\alpha_4)q_2,\\
r''_5:&x_5=-((q_1-1)p_1+(q_2-1)p_2-(\alpha_1+\alpha_5))p_1, \ y_5=\frac{1}{p_1},\\
&z_5=(q_2-1)p_1, \ w_5=\frac{p_2}{p_1},\\
r''_6:&x_6=-((q_1-t)p_1+(q_2-s)p_2-(\alpha_1+\alpha_6))p_1, \ y_6=\frac{1}{p_1},\\
&z_6=(q_2-s)p_1, \ w_6=\frac{p_2}{p_1}.
\end{split}
\end{align}
Then such a system coincides with the system \eqref{urauraS}.
\end{thm}
We remark that each transformation of each coordinate $r''_i \ (i=1,2,\dots,6),$ is birational and symplectic.

The difference between the conditions $r'_i$ and the conditions $r''_i$ is type of accessible singularities. More precisely, we will explain in the next section.

\begin{prop} 
On each coordinate $r''_j, \ j=1,2,\dots,6$, each of the  Hamiltonians $H''_{j1}$ and $H''_{j2}$ on $U_j \times B$ is expressed as a polynomial in $x_j,y_j,z_j,w_j$ and a rational function in $t$ and $s$, and satisfy the following conditions{\rm: \rm}
\begin{align}
\begin{split}
&dq_1 \wedge dp_1 +dq_2 \wedge dp_2 - dH''_1 \wedge dt- dH''_2 \wedge ds\\
&=dx_j \wedge dy_j +dz_j \wedge dw_j - dH''_{j1} \wedge dt- dH''_{j2} \wedge ds \quad (j=1,2,\dots,5),\\
&dq_1 \wedge dp_1 +dq_2 \wedge dp_2 - d(H''_1-p_1) \wedge dt- d(H''_2-p_2) \wedge ds\\
&=dx_6 \wedge dy_6 +dz_6 \wedge dw_6 - dH''_{61} \wedge dt- dH''_{62} \wedge ds.
\end{split}
\end{align}
\end{prop}

We also consider the inverse transformation of $S_1$.
\begin{Lemma}\label{th:ininokaoka}
The transformation
\begin{equation}
R_1:(q_1,p_1,q_2,p_2,t,s) \rightarrow \left(q_1,p_1-\frac{q_2p_2-\alpha_1-\alpha_3}{q_1},q_2q_1,\frac{p_2}{q_1},t,ts\right)
\end{equation}
satisfy the following relations\rm{:\rm}
\begin{equation}
R_1 \circ S_1=1, \quad S_1 \circ R_1=1.
\end{equation}
\end{Lemma}

Now, we also consider the following problem.
\begin{Problem}
For the system \eqref{urauraS}, can we find such a $S$-transformation differ-\\
ent from the transformation $R_1$?
\end{Problem}

To answer this, we find the following rational and symplectic transformation.
\begin{thm}\label{th:inininoka}
By using the rational and symplectic transformation $S_2$
\begin{equation}
S_2:(q_1,p_1,q_2,p_2,t,s) \rightarrow \left(q_1-\frac{q_2p_2-\alpha_1-\alpha_2-\alpha_3}{p_1},p_1,\frac{q_2}{p_1},p_2p_1\right),
\end{equation}
the system \eqref{urauraS} is transformed into a Hamiltonian system
\begin{equation}\label{uraurauraS}
\begin{aligned}
dq_1=\frac{\partial H_1'''}{\partial p_1}dt+\frac{\partial H_2'''}{\partial p_1}ds, \quad dp_1=-\frac{\partial H_1'''}{\partial q_1}dt-\frac{\partial H_2'''}{\partial q_1}ds,\\
dq_2=\frac{\partial H_1'''}{\partial p_2}dt+\frac{\partial H_2'''}{\partial p_2}ds, \quad dp_2=-\frac{\partial H_1'''}{\partial q_2}dt-\frac{\partial H_2'''}{\partial q_2}ds
\end{aligned}
\end{equation}
with polynomial Hamiltonians $H'''_1,H'''_2=\pi(H'''_1)$, where the transformation $\pi$ is explicitly given by
\begin{align}
\begin{split}
&\pi:(q_1,p_1,q_2,p_2,t,s;\alpha_1,\alpha_2, \dots ,\alpha_6)\\
&\rightarrow (-\frac{(q_1p_1-\alpha_2)p_1}{p_2},\frac{p_2}{p_1},q_2+\frac{q_1p_1-\alpha_1-\alpha_2-\alpha_3}{p_2},p_2,\\
&s,t;-\alpha_1-\alpha_2-\alpha_3-\alpha_4,\alpha_2,\alpha_4,\alpha_3,1-\alpha_6,1-\alpha_5 ).
\end{split}
\end{align}
\end{thm}
We note that each degree of each of Hamiltonians $H'''_1,H'''_2$ is 7 with respect to $q_1,p_1,q_2,p_2$.

\begin{thm}\label{holoURAURAURA}
Let us consider a polynomial Hamiltonian system with Hamiltonians $H_i \in {\Bbb C}(t,s)[q_1,p_1,q_2,p_2] \ (i=1,2)$. We assume that

$(A1)$ $deg(H_i)=7$ with respect to $q_1,p_1,q_2,p_2$.

$(A2)$ This system becomes again a polynomial Hamiltonian system in each coordinate $r'''_j, \ j=1,2,\dots,6${\rm : \rm}
\begin{align}
\begin{split}
r'''_1:&x_1=\frac{1}{q_1}, \ y_1=-(q_1p_1+q_2p_2-\alpha_1-\alpha_2)q_1, \ z_1=\frac{q_2}{q_1}, \ w_1=p_2q_1,\\
r'''_2:&x_2=-(q_1p_1-\alpha_2)p_1, \ y_2=\frac{1}{p_1}, \ z_2=q_2, \ w_2=p_2,\\
r'''_3:&x_3=\frac{1}{q_1},\ y_3=-(q_1p_1+q_2p_2-\alpha_1-\alpha_2-\alpha_3)q_1,\ z_3=\frac{q_2}{q_1},\ w_3=p_2q_1,\\
r'''_4:&x_4=q_1, \ y_4=p_1, \ z_4=\frac{1}{q_2}, \ w_4=-(p_2q_2+\alpha_4)q_2,\\
r'''_5:&x_5=-((q_1-1)p_1+2(q_2-\frac{1}{2p_1})p_2-(2\alpha_1+\alpha_2+\alpha_3+\alpha_5))p_1, \ y_5=\frac{1}{p_1},\\ 
&z_5=(q_2p_1-1)p_1, \ w_5=\frac{p_2}{p_1^2},\\
r'''_6:&x_6=-((q_1-t)p_1+2(q_2-\frac{s}{2p_1})p_2-(2\alpha_1+\alpha_2+\alpha_3+\alpha_6))p_1, \ y_6=\frac{1}{p_1},\\
&z_6=(q_2p_1-s)p_1, \ w_6=\frac{p_2}{p_1^2}.
\end{split}
\end{align}
Then such a system coincides with the system \eqref{uraurauraS}.
\end{thm}
We remark that each transformation of each coordinate $r'''_i \ (i=1,2,\dots,6),$ is birational and symplectic.

\begin{prop}
On each coordinate $r'''_j, \ j=1,2,\dots,6$, each of the  Hamiltonians $H'''_{j1}$ and $H'''_{j2}$ on $U_j \times B$ is expressed as a polynomial in $x_j,y_j,z_j,w_j$ and a rational function in $t$ and $s$, and satisfy the following conditions{\rm: \rm}
\begin{align}
\begin{split}
&dq_1 \wedge dp_1 +dq_2 \wedge dp_2 - dH'''_1 \wedge dt- dH'''_2 \wedge ds\\
&=dx_j \wedge dy_j +dz_j \wedge dw_j - dH'''_{j1} \wedge dt- dH'''_{j2} \wedge ds \quad (j=1,2,\dots,5),
\end{split}\\
\begin{split}
&dq_1 \wedge dp_1 +dq_2 \wedge dp_2 - d(H'''_1-p_1) \wedge dt- d(H'''_2-p_2/p_1) \wedge ds\\
&=dx_6 \wedge dy_6 +dz_6 \wedge dw_6 - dH'''_{61} \wedge dt- dH'''_{62} \wedge ds.
\end{split}
\end{align}
\end{prop}

We also consider the inverse transformation of $S_2$.
\begin{Lemma}\label{th:inininokaoka}
The transformation
\begin{equation}
R_2:(q_1,p_1,q_2,p_2) \rightarrow \left(q_1+\frac{q_2p_2-\alpha_1-\alpha_2-\alpha_3}{p_1},p_1,q_2p_1,\frac{p_2}{p_1}\right)
\end{equation}
satisfy the following relations\rm{:\rm}
\begin{equation}
R_2 \circ S_2=1, \quad S_2 \circ R_2=1.
\end{equation}
\end{Lemma}

\begin{Problem}
It is still an open question whether we classify such a $S$-transformat-\\
ion for each system.
\end{Problem}

\section{On some Hamiltonian structures of the system \eqref{uraS}}

In this section, we will explain the relation between the holomorphy conditions $r'_j$  with the accessible singularities of the system \eqref{uraS}. Let us take the compactification
\begin{equation*}
(q_1,p_1,q_2,p_2,t,s) \in {\Bbb C}^4 \times B_2 \  {\rm{to\rm}} \ ([z_0:z_1:z_2:z_3:z_4],t,s) \in {\Bbb P}^4 \times B_2
\end{equation*}
with the natural embedding
\begin{equation*}
(q_1,p_1,q_2,p_2)=(z_1/z_0,z_2/z_0,z_3/z_0,z_4/z_0).
\end{equation*}
Here $B_2={\Bbb C}^2-\{t(t-1)s(s-1)=0 \}$. Fixing the parameters $\alpha_i$, consider the product ${\Bbb P}^4 \times B_2$ and extend the regular vector field on ${\Bbb C}^4 \times B_2$ to a rational vector field $\tilde{v}$ on ${\Bbb P}^4 \times B_2$. It is easy to see that ${\Bbb P}^4$ is covered by five copies of ${\Bbb C}^4$ by gluing the transformations \eqref{coverofP4}. By the following lemma, we will show that each coordinate system $(x_i,y_i,z_i,w_i)$ $(i=1,2,\dots,6)$ can be obtained by successive blowing-up procedures of the accessible singularities in the boundary divisor ${\mathcal H} (\cong {\Bbb P}^3) \subset {\Bbb P}^4$.

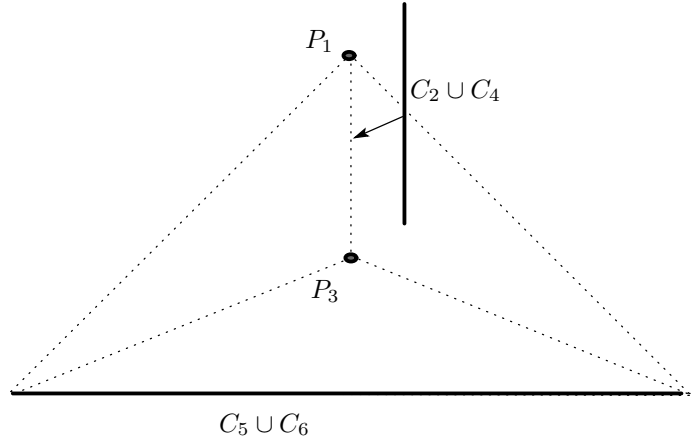
\begin{figure}[ht]
\unitlength 0.1in
\begin{picture}(35.50,20.70)(17.10,-23.60)
%
\special{pn 8}%
\special{pa 3490 550}%
\special{pa 1710 2320}%
\special{dt 0.045}%
\special{pa 1710 2320}%
\special{pa 1710 2320}%
\special{dt 0.045}%
%
\special{pn 8}%
\special{pa 3480 550}%
\special{pa 5260 2330}%
\special{dt 0.045}%
\special{pa 5260 2330}%
\special{pa 5260 2330}%
\special{dt 0.045}%
%
\special{pn 8}%
\special{pa 1730 2330}%
\special{pa 5260 2340}%
\special{dt 0.045}%
\special{pa 5260 2340}%
\special{pa 5259 2340}%
\special{dt 0.045}%
%
\special{pn 8}%
\special{pa 3490 560}%
\special{pa 3490 1630}%
\special{dt 0.045}%
\special{pa 3490 1630}%
\special{pa 3490 1629}%
\special{dt 0.045}%
%
\special{pn 8}%
\special{pa 1730 2330}%
\special{pa 3490 1620}%
\special{dt 0.045}%
\special{pa 3490 1620}%
\special{pa 3489 1620}%
\special{dt 0.045}%
%
\special{pn 8}%
\special{pa 3490 1610}%
\special{pa 5250 2330}%
\special{dt 0.045}%
\special{pa 5250 2330}%
\special{pa 5249 2330}%
\special{dt 0.045}%
%
\special{pn 20}%
\special{sh 0.600}%
\special{ar 3480 560 29 19  0.0000000 6.2831853}%
%
\special{pn 20}%
\special{sh 0.600}%
\special{ar 3490 1620 29 19  0.0000000 6.2831853}%
%
\special{pn 20}%
\special{pa 1720 2330}%
\special{pa 5220 2330}%
\special{fp}%
%
\special{pn 20}%
\special{pa 3770 290}%
\special{pa 3770 1440}%
\special{fp}%
%
\special{pn 8}%
\special{pa 3760 880}%
\special{pa 3510 990}%
\special{fp}%
\special{sh 1}%
\special{pa 3510 990}%
\special{pa 3579 981}%
\special{pa 3559 969}%
\special{pa 3563 945}%
\special{pa 3510 990}%
\special{fp}%
\put(32.5000,-5.2000){\makebox(0,0)[lb]{$P_1$}}%
\put(32.8000,-18.3000){\makebox(0,0)[lb]{$P_3$}}%
\put(38.0000,-7.8000){\makebox(0,0)[lb]{$C_2 \cup C_4$}}%
\put(28.0000,-25.3000){\makebox(0,0)[lb]{$C_5 \cup C_6$}}%
\end{picture}%
\label{fig:G2}
\caption{Accessible singularities of the system \eqref{uraS}}
\end{figure}

\begin{prop}\label{A-singularities of URA}
By resolving the following accessible singularities in the boundary divisor ${\mathcal H} (\cong {\Bbb P}^3) \subset {\Bbb P}^4$
\begin{equation}
  \left\{
  \begin{aligned}
   P_1 &=\{(X_1,Y_1,Z_1,W_1)|X_1=Y_1=Z_1=W_1=0 \},\\
   C_2 &=\{(X_1,Y_1,Z_1,W_1)|X_1=Y_1=W_1=0 \},\\
   P_3 &=\{(X_2,Y_2,Z_2,W_2)|X_2=Y_2=Z_2=W_2=0\},\\
   C_4 &=\{(X_2,Y_2,Z_2,W_2)|Y_2=Z_2=W_2=0\},\\
   C_5 &=\{(X_3,Y_3,Z_3,W_3)|X_3=Y_3=Z_3=0\},\\
   C_6 &=\{(X_4,Y_4,Z_4,W_4)|X_4=Z_4=W_4=0\},\\
   \end{aligned}
  \right. 
\end{equation}
we can obtain the coordinate systems $(x_i,y_i,z_i,w_i)$ $(i=1,2,\dots,6)$. Here $C_2 \cup C_4 \cong {\Bbb P}^1$, $C_5 \cup C_6 \cong {\Bbb P}^1$.
\end{prop}

\noindent
{\bf Proof of Proposition \ref{A-singularities of URA}.}
It is sufficient to show that we will resolve the accessible singularity $C_5 \cup C_6 \cong {\Bbb P}^1$. Other accessible singularities can be resolved as the same way in Proposition \ref{pro:asofG}.

{\bf Step 1:} We blow up along the curve $C_5 \cup C_6 \cong {\Bbb P}^1${\rm : \rm}
$$
{x_{1}^{(1)}}=\frac{X_3}{Y_3} \;, \;\;\; {y_{1}^{(1)}}=Y_3 \;, \;\;\; {z_{1}^{(1)}}=\frac{Z_3}{Y_3} \;, \;\;\; {w_{1}^{(1)}}=W_3,
$$
$$
{x_{2}^{(1)}}=\frac{X_4}{W_4} \;, \;\;\; {y_{2}^{(1)}}=Y_4 \;, \;\;\; {z_{2}^{(1)}}=\frac{Z_4}{W_4} \;, \;\;\; {w_{2}^{(1)}}=W_4.
$$

It is easy to see that there are three accessible singularities
\begin{align*}
L_1=&\{(x_{1}^{(1)},y_{1}^{(1)},z_{1}^{(1)},w_{1}^{(1)})|x_{1}^{(1)}=y_{1}^{(1)}=z_{1}^{(1)}=0 \} \cup\\
&\{(x_{2}^{(1)},y_{2}^{(1)},z_{2}^{(1)},w_{2}^{(1)})|x_{2}^{(1)}=z_{2}^{(1)}=w_{2}^{(1)}=0 \} \cong {\Bbb P}^1,\\
L_2=&\{(x_{1}^{(1)},y_{1}^{(1)},z_{1}^{(1)},w_{1}^{(1)})| x_{1}^{(1)}-1=y_{1}^{(1)}=z_{1}^{(1)}-1=0 \} \cup\\
&\{(x_{2}^{(1)},y_{2}^{(1)},z_{2}^{(1)},w_{2}^{(1)})|x_{2}^{(1)}-1=z_{2}^{(1)}-1=w_{2}^{(1)}=0 \} \cong {\Bbb P}^1,\\
L_3=&\{(x_{1}^{(1)},y_{1}^{(1)},z_{1}^{(1)},w_{1}^{(1)})|x_{1}^{(1)}-t=y_{1}^{(1)}=z_{1}^{(1)}-s=0 \} \cup\\
&\{(x_{2}^{(1)},y_{2}^{(1)},z_{2}^{(1)},w_{2}^{(1)})|x_{2}^{(1)}-t=z_{2}^{(1)}-s=w_{2}^{(1)}=0 \} \cong {\Bbb P}^1.
\end{align*}

{\bf Step 2:} We blow up along the curve $L_1${\rm : \rm}
\begin{align*}
&{x_{1}^{(2)}}=\frac{x_{1}^{(1)}}{y_{1}^{(1)}} \;, \;\;\; {y_{1}^{(2)}}=y_{1}^{(1)} \;, \;\;\; {z_{1}^{(2)}}=\frac{z_{1}^{(1)}}{y_{1}^{(1)}} \;, \;\;\; {w_{1}^{(2)}}=w_{1}^{(1)},\\
&{x_{2}^{(2)}}=\frac{x_{2}^{(1)}}{w_{2}^{(1)}} \;, \;\;\; {y_{2}^{(2)}}=y_{2}^{(1)} \;, \;\;\; {z_{2}^{(2)}}=\frac{z_{2}^{(1)}}{w_{2}^{(1)}} \;, \;\;\; {w_{2}^{(2)}}=w_{2}^{(1)}.
\end{align*}

{\bf Step 3:} We blow up along the surface
\begin{align*}
S_1=&\{(x_{1}^{(2)},y_{1}^{(2)},z_{1}^{(2)},w_{1}^{(2)})|x_{1}^{(2)}+z_{1}^{(2)}w_{1}^{(2)}-(\alpha_1+\alpha_3)=y_{1}^{(2)}=0 \} \cup\\
&\{(x_{2}^{(2)},y_{2}^{(2)},z_{2}^{(2)},w_{2}^{(2)})|z_{2}^{(2)}+x_{2}^{(2)}y_{2}^{(2)}-(\alpha_1+\alpha_3)=w_{2}^{(2)}=0 \}{\rm : \rm}
\end{align*}
\begin{align*}
&{x_{1}^{(3)}}=\frac{x_{1}^{(2)}+z_{1}^{(2)}w_{1}^{(2)}-(\alpha_1+\alpha_3)}{y_{1}^{(2)}} \;, \;\;\; {y_{1}^{(3)}}=y_{1}^{(2)} \;, \;\;\; {z_{1}^{(3)}}=z_{1}^{(2)} \;, \;\;\; {w_{1}^{(3)}}=w_{1}^{(2)},\\
&{x_{2}^{(3)}}=x_{2}^{(2)} \;, \;\;\; {y_{2}^{(3)}}=y_{2}^{(2)} \;, \;\;\; {z_{2}^{(3)}}=\frac{z_{2}^{(2)}+x_{2}^{(2)}y_{2}^{(2)}-(\alpha_1+\alpha_3)}{w_{2}^{(2)}} \;, \;\;\; {w_{2}^{(3)}}=w_{2}^{(2)}.
\end{align*}
We have resolved the accessible singularity $L_1$.

For the remaining accessible singularities, the proof is similar.

The proof has thus been completed.   \qed

By the same way, we will explain the relation between the holomorphy conditions $r''_j$ with the following accessible singularities of the system \eqref{urauraS}.

\begin{figure}
\unitlength 0.1in
\begin{picture}(35.50,20.10)(17.10,-23.60)
%
\special{pn 8}%
\special{pa 3490 550}%
\special{pa 1710 2320}%
\special{dt 0.045}%
\special{pa 1710 2320}%
\special{pa 1710 2320}%
\special{dt 0.045}%
%
\special{pn 20}%
\special{pa 3480 550}%
\special{pa 5260 2330}%
\special{fp}%
%
\special{pn 8}%
\special{pa 1730 2330}%
\special{pa 5260 2340}%
\special{dt 0.045}%
\special{pa 5260 2340}%
\special{pa 5259 2340}%
\special{dt 0.045}%
%
\special{pn 8}%
\special{pa 3490 560}%
\special{pa 3490 1630}%
\special{dt 0.045}%
\special{pa 3490 1630}%
\special{pa 3490 1629}%
\special{dt 0.045}%
%
\special{pn 20}%
\special{pa 1730 2330}%
\special{pa 3490 1620}%
\special{fp}%
%
\special{pn 8}%
\special{pa 3490 1610}%
\special{pa 5250 2330}%
\special{dt 0.045}%
\special{pa 5250 2330}%
\special{pa 5249 2330}%
\special{dt 0.045}%
%
\special{pn 20}%
\special{sh 0.600}%
\special{ar 3480 560 29 19  0.0000000 6.2831853}%
%
\special{pn 20}%
\special{sh 0.600}%
\special{ar 3490 1620 29 19  0.0000000 6.2831853}%
%
\special{pn 20}%
\special{pa 1720 2330}%
\special{pa 5220 2330}%
\special{fp}%
\put(32.5000,-5.2000){\makebox(0,0)[lb]{$P_3$}}%
\put(32.8000,-18.3000){\makebox(0,0)[lb]{$P_4$}}%
\put(37.3000,-14.4000){\makebox(0,0)[lb]{$C_5 \cup C_6$}}%
\put(28.0000,-25.3000){\makebox(0,0)[lb]{$C_1 \cup C_2$}}%
\put(21.9000,-18.7000){\makebox(0,0)[lb]{$C_7 \cup C_8$}}%
\end{picture}%
\label{fig:G3}
\caption{Accessible singularities of the system \eqref{urauraS}}
\end{figure}
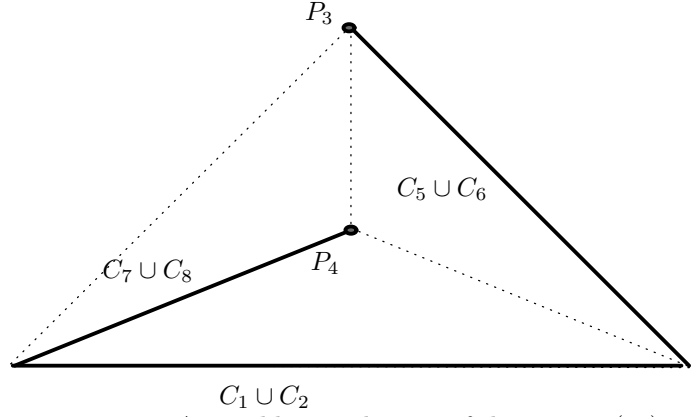

\begin{prop}\label{Acceuraura}
By resolving the following accessible singularities in the boundary divisor ${\mathcal H} (\cong {\Bbb P}^3) \subset {\Bbb P}^4$
\begin{equation}
  \left\{
  \begin{aligned}
   C_1 &=\{(X_3,Y_3,Z_3,W_3)|X_3=Y_3=Z_3=0\},\\
   C_2 &=\{(X_4,Y_4,Z_4,W_4)|X_4=Z_4=W_4=0\},\\
   P_3 &=\{(X_1,Y_1,Z_1,W_1)|X_1=Y_1=Z_1=W_1=0 \},\\
   P_4 &=\{(X_2,Y_2,Z_2,W_2)|X_2=Y_2=Z_2=W_2=0\},\\
   C_5 &=\{(X_1,Y_1,Z_1,W_1)|X_1=Y_1=Z_1=0\},\\
   C_6 &=\{(X_4,Y_4,Z_4,W_4)|Y_4=Z_4=W_4=0\},\\
   C_7 &=\{(X_2,Y_2,Z_2,W_2)|X_2=Z_2=W_2=0 \},\\
   C_8 &=\{(X_3,Y_3,Z_3,W_3)|X_3=Y_3=W_3=0 \},\\
   \end{aligned}
  \right. 
\end{equation}
we can obtain the coordinates $r''_i \ (i=1,2,\dots,6)$. Here $C_1 \cup C_2 \cong {\Bbb P}^1$, $C_5 \cup C_6 \cong {\Bbb P}^1$ and $C_7 \cup C_8 \cong {\Bbb P}^1$.
\end{prop}
We remark that the relations between the accessible singularities given by Proposition \ref{Acceuraura} and the coordinates $r''_j \ (j=1,2,\dots,6)$ are given as follows:
\begin{center}
\begin{tabular}{|c|c|c|c|c|c|} \hline 
$C_1 \cup C_2$ &$C_1 \cup C_2$  & $P_3$ & $P_4$ & $C_5 \cup C_6$& $C_7 \cup C_8$  \\ \hline 
$r''_5$ & $r''_6$ & $r''_3$ & $r''_4$ & $r''_2$ & $r''_1$  \\ \hline 
\end{tabular}
\end{center}

\section{Invariant cycles of the system \eqref{1}}

\begin{figure}
\unitlength 0.1in
\begin{picture}(37.19,34.02)(19.96,-34.76)
%
\special{pn 20}%
\special{ar 4054 893 293 193  0.0000000 6.2831853}%
%
\special{pn 20}%
\special{ar 5436 234 231 160  0.0000000 6.2831853}%
%
\special{pn 20}%
\special{ar 5484 1565 231 160  0.0000000 6.2831853}%
%
\special{pn 8}%
\special{pa 4217 718}%
\special{pa 5206 248}%
\special{fp}%
%
\special{pn 8}%
\special{pa 4217 1063}%
\special{pa 5244 1515}%
\special{fp}%
%
\special{pn 8}%
\special{pa 3766 907}%
\special{pa 4342 907}%
\special{fp}%
\put(20.7000,-7.0000){\makebox(0,0)[lb]{$q_1+q_2-1$}}%
\put(20.2000,-12.2000){\makebox(0,0)[lb]{$q_1+tq_2/s-t$}}%
\put(39.1000,-8.8100){\makebox(0,0)[lb]{$p_1$}}%
\put(39.2000,-10.2400){\makebox(0,0)[lb]{$p_2$}}%
\put(53.0200,-2.9200){\makebox(0,0)[lb]{$q_1$}}%
\put(53.4000,-16.2700){\makebox(0,0)[lb]{$q_2$}}%
%
\special{pn 20}%
\special{ar 2441 1175 426 162  0.0000000 6.2831853}%
%
\special{pn 20}%
\special{ar 2422 661 426 161  0.0000000 6.2831853}%
%
\special{pn 20}%
\special{ar 2677 2662 267 245  0.0000000 6.2831853}%
%
\special{pn 20}%
\special{ar 3677 2668 267 244  0.0000000 6.2831853}%
%
\special{pn 20}%
\special{ar 4677 2675 267 243  0.0000000 6.2831853}%
%
\special{pn 20}%
\special{ar 3149 1966 192 146  0.0000000 6.2831853}%
%
\special{pn 20}%
\special{ar 3165 3330 192 146  0.0000000 6.2831853}%
%
\special{pn 8}%
\special{pa 3037 2086}%
\special{pa 2741 2423}%
\special{fp}%
%
\special{pn 8}%
\special{pa 3269 2086}%
\special{pa 3589 2429}%
\special{fp}%
%
\special{pn 8}%
\special{pa 3341 1994}%
\special{pa 4533 2451}%
\special{fp}%
%
\special{pn 8}%
\special{pa 2773 2896}%
\special{pa 3013 3209}%
\special{fp}%
%
\special{pn 8}%
\special{pa 4509 2873}%
\special{pa 3349 3278}%
\special{fp}%
%
\special{pn 8}%
\special{pa 3533 2873}%
\special{pa 3285 3204}%
\special{fp}%
%
\special{pn 8}%
\special{pa 2421 2646}%
\special{pa 2925 2646}%
\special{fp}%
%
\special{pn 8}%
\special{pa 3413 2657}%
\special{pa 3925 2657}%
\special{fp}%
%
\special{pn 8}%
\special{pa 4421 2668}%
\special{pa 4933 2668}%
\special{fp}%
\put(30.6900,-20.0700){\makebox(0,0)[lb]{$p_1$}}%
\put(30.9300,-34.0400){\makebox(0,0)[lb]{$p_2$}}%
\put(25.4100,-25.9400){\makebox(0,0)[lb]{$q_1$}}%
\put(25.4100,-27.8200){\makebox(0,0)[lb]{$q_2$}}%
\put(35.0900,-27.7600){\makebox(0,0)[lb]{$q_2-s$}}%
\put(44.9300,-26.1100){\makebox(0,0)[lb]{$q_1-1$}}%
\put(45.0100,-27.7600){\makebox(0,0)[lb]{$q_2-1$}}%
%
\special{pn 8}%
\special{pa 2844 674}%
\special{pa 3747 913}%
\special{fp}%
%
\special{pn 8}%
\special{pa 2864 1195}%
\special{pa 3737 907}%
\special{fp}%
\put(35.0900,-26.0500){\makebox(0,0)[lb]{$q_1-t$}}%
\end{picture}%
\label{UraDynkin}
\caption{The symbol in each circle denotes the invariant cycle for each system.}
\end{figure}
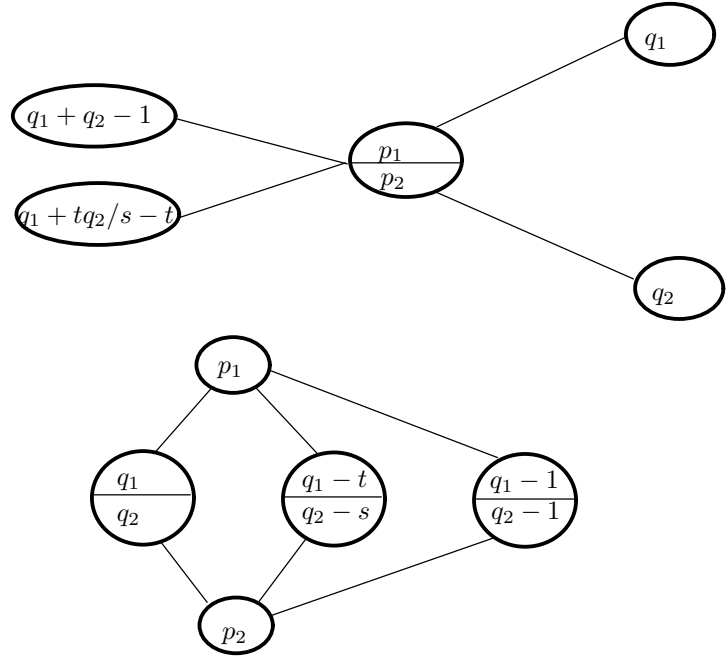

\begin{figure}
\unitlength 0.1in
\begin{picture}(49.09,28.80)(15.20,-31.80)
%
\special{pn 8}%
\special{pa 1917 514}%
\special{pa 2449 1494}%
\special{dt 0.045}%
\special{pa 2449 1494}%
\special{pa 2449 1493}%
\special{dt 0.045}%
\special{pa 2449 1494}%
\special{pa 3032 1701}%
\special{dt 0.045}%
\special{pa 3032 1701}%
\special{pa 3031 1701}%
\special{dt 0.045}%
\special{pa 3032 1701}%
\special{pa 3933 1696}%
\special{dt 0.045}%
\special{pa 3933 1696}%
\special{pa 3932 1696}%
\special{dt 0.045}%
\special{pa 1917 514}%
\special{pa 3928 1690}%
\special{dt 0.045}%
\special{pa 3928 1690}%
\special{pa 3927 1690}%
\special{dt 0.045}%
%
\special{pn 8}%
\special{pa 1923 520}%
\special{pa 3015 1696}%
\special{dt 0.045}%
\special{pa 3015 1696}%
\special{pa 3015 1696}%
\special{dt 0.045}%
%
\special{pn 13}%
\special{pa 1923 514}%
\special{pa 1928 2564}%
\special{fp}%
\special{pa 1928 2564}%
\special{pa 3922 1707}%
\special{fp}%
%
\special{pn 13}%
\special{pa 1928 2547}%
\special{pa 2444 1494}%
\special{fp}%
%
\special{pn 13}%
\special{pa 1928 2552}%
\special{pa 3015 1701}%
\special{fp}%
%
\special{pn 8}%
\special{pa 2444 1488}%
\special{pa 3916 1690}%
\special{dt 0.045}%
\special{pa 3916 1690}%
\special{pa 3915 1690}%
\special{dt 0.045}%
%
\special{pn 8}%
\special{pa 2281 1197}%
\special{pa 1928 1550}%
\special{fp}%
\special{pa 2292 1220}%
\special{pa 1928 1584}%
\special{fp}%
\special{pa 2304 1242}%
\special{pa 1928 1617}%
\special{fp}%
\special{pa 2320 1259}%
\special{pa 1928 1651}%
\special{fp}%
\special{pa 2332 1281}%
\special{pa 1928 1684}%
\special{fp}%
\special{pa 2343 1304}%
\special{pa 1928 1718}%
\special{fp}%
\special{pa 2354 1326}%
\special{pa 1928 1752}%
\special{fp}%
\special{pa 2365 1348}%
\special{pa 1928 1785}%
\special{fp}%
\special{pa 2376 1371}%
\special{pa 1928 1819}%
\special{fp}%
\special{pa 2388 1393}%
\special{pa 1928 1852}%
\special{fp}%
\special{pa 2399 1416}%
\special{pa 1928 1886}%
\special{fp}%
\special{pa 2410 1438}%
\special{pa 1928 1920}%
\special{fp}%
\special{pa 2427 1455}%
\special{pa 1928 1953}%
\special{fp}%
\special{pa 2438 1477}%
\special{pa 1928 1987}%
\special{fp}%
\special{pa 2427 1522}%
\special{pa 1928 2020}%
\special{fp}%
\special{pa 2393 1589}%
\special{pa 1928 2054}%
\special{fp}%
\special{pa 2360 1656}%
\special{pa 1928 2088}%
\special{fp}%
\special{pa 2326 1724}%
\special{pa 1928 2121}%
\special{fp}%
\special{pa 2298 1785}%
\special{pa 1928 2155}%
\special{fp}%
\special{pa 2264 1852}%
\special{pa 1928 2188}%
\special{fp}%
\special{pa 2231 1920}%
\special{pa 1928 2222}%
\special{fp}%
\special{pa 2197 1987}%
\special{pa 1928 2256}%
\special{fp}%
\special{pa 2169 2048}%
\special{pa 1928 2289}%
\special{fp}%
\special{pa 2136 2116}%
\special{pa 1928 2323}%
\special{fp}%
\special{pa 2102 2183}%
\special{pa 1928 2356}%
\special{fp}%
\special{pa 2068 2250}%
\special{pa 1928 2390}%
\special{fp}%
\special{pa 2040 2312}%
\special{pa 1928 2424}%
\special{fp}%
\special{pa 2007 2379}%
\special{pa 1928 2457}%
\special{fp}%
\special{pa 1973 2446}%
\special{pa 1928 2491}%
\special{fp}%
\special{pa 2270 1175}%
\special{pa 1928 1516}%
\special{fp}%
%
\special{pn 8}%
\special{pa 2259 1152}%
\special{pa 1928 1483}%
\special{fp}%
\special{pa 2248 1130}%
\special{pa 1928 1449}%
\special{fp}%
\special{pa 2236 1108}%
\special{pa 1928 1416}%
\special{fp}%
\special{pa 2225 1085}%
\special{pa 1928 1382}%
\special{fp}%
\special{pa 2214 1063}%
\special{pa 1928 1348}%
\special{fp}%
\special{pa 2197 1046}%
\special{pa 1928 1315}%
\special{fp}%
\special{pa 2186 1024}%
\special{pa 1928 1281}%
\special{fp}%
\special{pa 2175 1001}%
\special{pa 1928 1248}%
\special{fp}%
\special{pa 2164 979}%
\special{pa 1928 1214}%
\special{fp}%
\special{pa 2152 956}%
\special{pa 1928 1180}%
\special{fp}%
\special{pa 2141 934}%
\special{pa 1928 1147}%
\special{fp}%
\special{pa 2130 912}%
\special{pa 1928 1113}%
\special{fp}%
\special{pa 2119 889}%
\special{pa 1928 1080}%
\special{fp}%
\special{pa 2108 867}%
\special{pa 1928 1046}%
\special{fp}%
\special{pa 2091 850}%
\special{pa 1928 1012}%
\special{fp}%
\special{pa 2080 828}%
\special{pa 1928 979}%
\special{fp}%
\special{pa 2068 805}%
\special{pa 1928 945}%
\special{fp}%
\special{pa 2057 783}%
\special{pa 1928 912}%
\special{fp}%
\special{pa 2046 760}%
\special{pa 1928 878}%
\special{fp}%
\special{pa 2035 738}%
\special{pa 1928 844}%
\special{fp}%
\special{pa 2024 716}%
\special{pa 1928 811}%
\special{fp}%
\special{pa 2012 693}%
\special{pa 1928 777}%
\special{fp}%
\special{pa 2001 671}%
\special{pa 1928 744}%
\special{fp}%
\special{pa 1984 654}%
\special{pa 1928 710}%
\special{fp}%
\special{pa 1973 632}%
\special{pa 1928 676}%
\special{fp}%
\special{pa 1962 609}%
\special{pa 1928 643}%
\special{fp}%
\special{pa 1951 587}%
\special{pa 1928 609}%
\special{fp}%
\special{pa 1940 564}%
\special{pa 1928 576}%
\special{fp}%
%
\special{pn 20}%
\special{pa 1928 2558}%
\special{pa 1928 2328}%
\special{fp}%
\special{sh 1}%
\special{pa 1928 2328}%
\special{pa 1908 2395}%
\special{pa 1928 2381}%
\special{pa 1948 2395}%
\special{pa 1928 2328}%
\special{fp}%
%
\special{pn 20}%
\special{pa 1923 2547}%
\special{pa 2052 2306}%
\special{fp}%
\special{sh 1}%
\special{pa 2052 2306}%
\special{pa 2003 2355}%
\special{pa 2027 2353}%
\special{pa 2038 2374}%
\special{pa 2052 2306}%
\special{fp}%
%
\special{pn 20}%
\special{pa 1923 2547}%
\special{pa 2180 2351}%
\special{fp}%
\special{sh 1}%
\special{pa 2180 2351}%
\special{pa 2115 2376}%
\special{pa 2138 2383}%
\special{pa 2139 2407}%
\special{pa 2180 2351}%
\special{fp}%
%
\special{pn 20}%
\special{pa 1934 2558}%
\special{pa 2192 2452}%
\special{fp}%
\special{sh 1}%
\special{pa 2192 2452}%
\special{pa 2123 2459}%
\special{pa 2143 2472}%
\special{pa 2138 2496}%
\special{pa 2192 2452}%
\special{fp}%
\put(18.5000,-22.9500){\makebox(0,0)[lb]{$q_1$}}%
\put(20.2400,-22.9500){\makebox(0,0)[lb]{$q_2$}}%
\put(21.5800,-23.4500){\makebox(0,0)[lb]{$p_1$}}%
\put(21.8000,-24.7400){\makebox(0,0)[lb]{$p_2$}}%
%
\special{pn 8}%
\special{pa 4413 490}%
\special{pa 4945 1470}%
\special{dt 0.045}%
\special{pa 4945 1470}%
\special{pa 4945 1469}%
\special{dt 0.045}%
\special{pa 4945 1470}%
\special{pa 5528 1677}%
\special{dt 0.045}%
\special{pa 5528 1677}%
\special{pa 5527 1677}%
\special{dt 0.045}%
\special{pa 5528 1677}%
\special{pa 6429 1672}%
\special{dt 0.045}%
\special{pa 6429 1672}%
\special{pa 6428 1672}%
\special{dt 0.045}%
\special{pa 4413 490}%
\special{pa 6424 1666}%
\special{dt 0.045}%
\special{pa 6424 1666}%
\special{pa 6423 1666}%
\special{dt 0.045}%
%
\special{pn 8}%
\special{pa 4419 496}%
\special{pa 5511 1672}%
\special{dt 0.045}%
\special{pa 5511 1672}%
\special{pa 5511 1672}%
\special{dt 0.045}%
%
\special{pn 13}%
\special{pa 4419 490}%
\special{pa 4424 2540}%
\special{fp}%
\special{pa 4424 2540}%
\special{pa 6418 1683}%
\special{fp}%
%
\special{pn 13}%
\special{pa 4424 2523}%
\special{pa 4940 1470}%
\special{fp}%
%
\special{pn 13}%
\special{pa 4424 2528}%
\special{pa 5511 1677}%
\special{fp}%
%
\special{pn 8}%
\special{pa 4940 1464}%
\special{pa 6412 1666}%
\special{dt 0.045}%
\special{pa 6412 1666}%
\special{pa 6411 1666}%
\special{dt 0.045}%
%
\special{pn 4}%
\special{pa 5785 1677}%
\special{pa 5309 2153}%
\special{fp}%
\special{pa 5752 1677}%
\special{pa 5248 2181}%
\special{fp}%
\special{pa 5718 1677}%
\special{pa 5192 2204}%
\special{fp}%
\special{pa 5684 1677}%
\special{pa 5130 2232}%
\special{fp}%
\special{pa 5651 1677}%
\special{pa 5074 2254}%
\special{fp}%
\special{pa 5617 1677}%
\special{pa 5012 2282}%
\special{fp}%
\special{pa 5584 1677}%
\special{pa 4956 2304}%
\special{fp}%
\special{pa 5550 1677}%
\special{pa 4895 2332}%
\special{fp}%
\special{pa 5505 1688}%
\special{pa 4839 2355}%
\special{fp}%
\special{pa 5348 1812}%
\special{pa 4777 2383}%
\special{fp}%
\special{pa 5192 1935}%
\special{pa 4716 2411}%
\special{fp}%
\special{pa 5040 2052}%
\special{pa 4660 2433}%
\special{fp}%
\special{pa 4884 2176}%
\special{pa 4598 2461}%
\special{fp}%
\special{pa 4727 2299}%
\special{pa 4542 2484}%
\special{fp}%
\special{pa 4570 2422}%
\special{pa 4480 2512}%
\special{fp}%
\special{pa 5819 1677}%
\special{pa 5365 2131}%
\special{fp}%
\special{pa 5852 1677}%
\special{pa 5427 2103}%
\special{fp}%
\special{pa 5886 1677}%
\special{pa 5483 2080}%
\special{fp}%
\special{pa 5920 1677}%
\special{pa 5544 2052}%
\special{fp}%
\special{pa 5953 1677}%
\special{pa 5600 2030}%
\special{fp}%
\special{pa 5987 1677}%
\special{pa 5662 2002}%
\special{fp}%
\special{pa 6020 1677}%
\special{pa 5718 1980}%
\special{fp}%
\special{pa 6054 1677}%
\special{pa 5780 1952}%
\special{fp}%
\special{pa 6088 1677}%
\special{pa 5836 1929}%
\special{fp}%
\special{pa 6121 1677}%
\special{pa 5897 1901}%
\special{fp}%
\special{pa 6155 1677}%
\special{pa 5953 1879}%
\special{fp}%
\special{pa 6188 1677}%
\special{pa 6015 1851}%
\special{fp}%
\special{pa 6222 1677}%
\special{pa 6071 1828}%
\special{fp}%
\special{pa 6256 1677}%
\special{pa 6132 1800}%
\special{fp}%
\special{pa 6289 1677}%
\special{pa 6194 1772}%
\special{fp}%
%
\special{pn 4}%
\special{pa 6323 1677}%
\special{pa 6250 1750}%
\special{fp}%
\special{pa 6356 1677}%
\special{pa 6312 1722}%
\special{fp}%
%
\special{pn 20}%
\special{pa 4424 2534}%
\special{pa 4424 2304}%
\special{fp}%
\special{sh 1}%
\special{pa 4424 2304}%
\special{pa 4404 2371}%
\special{pa 4424 2357}%
\special{pa 4444 2371}%
\special{pa 4424 2304}%
\special{fp}%
%
\special{pn 20}%
\special{pa 4419 2523}%
\special{pa 4548 2282}%
\special{fp}%
\special{sh 1}%
\special{pa 4548 2282}%
\special{pa 4499 2331}%
\special{pa 4523 2329}%
\special{pa 4534 2350}%
\special{pa 4548 2282}%
\special{fp}%
%
\special{pn 20}%
\special{pa 4419 2523}%
\special{pa 4676 2327}%
\special{fp}%
\special{sh 1}%
\special{pa 4676 2327}%
\special{pa 4611 2352}%
\special{pa 4634 2359}%
\special{pa 4635 2383}%
\special{pa 4676 2327}%
\special{fp}%
%
\special{pn 20}%
\special{pa 4430 2534}%
\special{pa 4688 2428}%
\special{fp}%
\special{sh 1}%
\special{pa 4688 2428}%
\special{pa 4619 2435}%
\special{pa 4639 2448}%
\special{pa 4634 2472}%
\special{pa 4688 2428}%
\special{fp}%
\put(43.4600,-22.7100){\makebox(0,0)[lb]{$q_1$}}%
\put(45.2000,-22.7100){\makebox(0,0)[lb]{$q_2$}}%
\put(46.5400,-23.2100){\makebox(0,0)[lb]{$p_1$}}%
\put(46.7600,-24.5000){\makebox(0,0)[lb]{$p_2$}}%
\put(22.4000,-29.6000){\makebox(0,0)[lb]{Type A}}%
\put(47.5000,-29.4000){\makebox(0,0)[lb]{Type B}}%
%
\special{pn 20}%
\special{sh 0.600}%
\special{ar 1930 510 15 23  0.0000000 6.2831853}%
%
\special{pn 20}%
\special{sh 0.600}%
\special{ar 2430 1490 15 23  0.0000000 6.2831853}%
%
\special{pn 20}%
\special{sh 0.600}%
\special{ar 3010 1690 15 23  0.0000000 6.2831853}%
%
\special{pn 20}%
\special{sh 0.600}%
\special{ar 3900 1700 15 23  0.0000000 6.2831853}%
\put(17.5000,-4.7000){\makebox(0,0)[lb]{$P_1$}}%
\put(23.3000,-14.3000){\makebox(0,0)[lb]{$P_2$}}%
\put(28.8000,-16.5000){\makebox(0,0)[lb]{$P_3$}}%
\put(37.4000,-16.6000){\makebox(0,0)[lb]{$P_4$}}%
%
\put(15.2000,-31.8000){\makebox(0,0)[lb]{}}%
%
\special{pn 20}%
\special{sh 0.600}%
\special{ar 1940 2550 15 23  0.0000000 6.2831853}%
\put(17.8000,-27.1000){\makebox(0,0)[lb]{$P_0$}}%
%
\special{pn 20}%
\special{pa 1940 520}%
\special{pa 2440 1490}%
\special{fp}%
%
\special{pn 20}%
\special{pa 5510 1680}%
\special{pa 6380 1680}%
\special{fp}%
\end{picture}%
\label{invariantURA}
\caption{This figure denotes the four-dimensional projective space ${\Bbb P}^4={\Bbb C}^4 \sqcup {\Bbb P}^3$. ${\Bbb P}^4$ is covered by five open sets ${\Bbb C}^4$ around the points $P_i \ (i=0,1,\dots,4)$. We also remark that the figure spanned by the points $P_i \ (i=1,2,\dots,4)$ denotes the three-dimensional projective space ${\Bbb P}^3={\Bbb C}^3 \sqcup {\Bbb P}^2$, and spanned by the points $P_i \ (i=0,1,2)$ (resp. $P_i \ (i=0,3,4)$) denotes the two-dimensional projective space ${\Bbb P}^2={\Bbb C}^2 \sqcup {\Bbb P}^1$.
}
\end{figure}
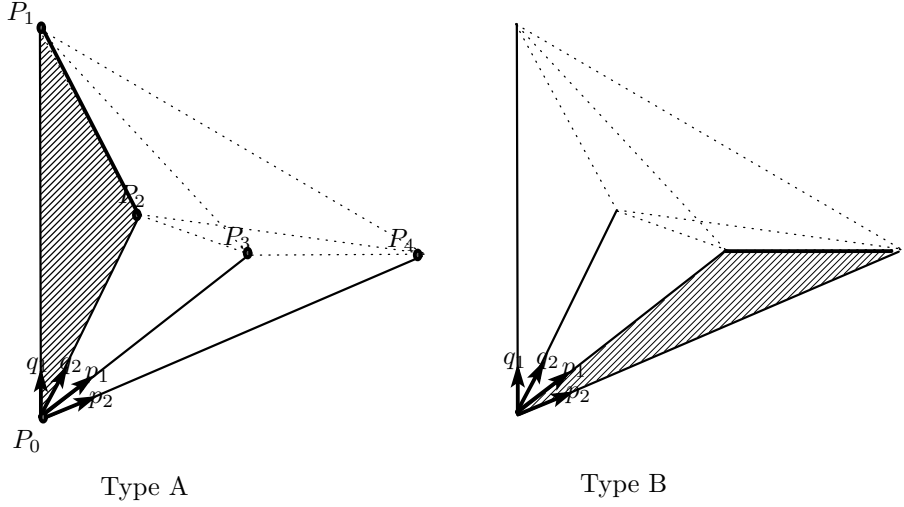

\begin{center}
\begin{tabular}{|c|c|c|} \hline
codimension & invariant cycle & parameter's relation \\ \hline
2 & $f_1^{(1)}:=p_1, \quad f_1^{(2)}:=p_2$ & $\alpha_1=0$  \\ \hline
2 & $f_2^{(1)}:=p_1, \quad f_2^{(2)}:=p_2$ & $\alpha_1=-\alpha_2$  \\ \hline
1 & $f_3:=q_1$ & $\alpha_3=0$  \\ \hline
1 & $f_4:=q_2$ & $\alpha_4=0$  \\ \hline
1 & $f_5:=q_1+q_2-1$ & $\alpha_5=0$   \\ \hline
1 & $f_6:=q_1+\frac{tq_2}{s}-t$ & $\alpha_6=0$   \\ \hline
\end{tabular}
\end{center}
The list must be read as follows:

\noindent
Setting $\alpha_1=0$, then we see that the system \eqref{1} admits a particular solution expressed in terms of Appell's hypergeometric function:
\begin{align}
\begin{split}
&p_1=p_2=0,\\
dq_1=&\{-\frac{\alpha_3(q_1-1)(q_1-t)+(\alpha_4+\alpha_6-1)q_1(q_1-1)+\alpha_5q_1(q_1-t)}{t(t-1)}\\
&-\frac{\alpha_4 sq_1}{t(t-s)}+\frac{\alpha_3 q_2}{t-s}\}dt+\{\frac{\alpha_2 q_1q_2}{s(s-1)}+\frac{\alpha_3 tq_2}{s(s-t)}-\frac{\alpha_4(t-1)q_1}{(s-1)(s-t)}\}ds,\\
dq_2=&\{\frac{\alpha_2 q_1q_2}{t(t-1)}+\frac{\alpha_4 sq_1}{t(t-s)}-\frac{\alpha_3(s-1)q_2}{(t-1)(t-s)}\}dt\\
&+\{-\frac{\alpha_4(q_2-1)(q_2-s)+(\alpha_3+\alpha_6-1)q_2(q_2-1)+\alpha_5q_2(q_2-s)}{s(s-1)}\\
&-\frac{\alpha_3 tq_2}{s(s-t)}+\frac{\alpha_4 q_1}{s-t}\}ds.
\end{split}
\end{align}
This system is invariant under the transformation $\pi$:
$$
\pi:(q_1,q_2,t,s;\alpha_2,\dots,\alpha_6) \rightarrow (q_2,q_1,s,t;\alpha_2,\alpha_4,\alpha_3,\alpha_5,\alpha_6).
$$
This case corresponds to the slant part of Type A in Figure 6.

\noindent
And, setting $\alpha_3=0$, then the system \eqref{1} admits a particular solution $q_1=0$. Moreover $(q_2,p_2)$ satisfy the sixth Painlev\'e system $H_{VI}$. And $p_1$ satisfies Riccati type equations whose coefficients are polynomials $(p_1,p_2)$, and so on.

\section{Invariant cycles of the system \eqref{uraS}}

\begin{center}
\begin{tabular}{|c|c|c|} \hline
codimension & invariant cycle & parameter's relation \\ \hline
1 & $f_1:=p_1$ & $\alpha_2=0$  \\ \hline
1 & $f_3:=p_2$ & $\alpha_4=0$  \\ \hline
2 & $f_4^{(1)}:=q_1, \ f_4^{(2)}:=q_2$ & $\alpha_3=-\alpha_1$  \\ \hline
2 & $f_5^{(1)}:=q_1-1, \ f_5^{(2)}:=q_2-1$ & $\alpha_5=-\alpha_1$   \\ \hline
2 & $f_6^{(1)}:=q_1-t, \ f_6^{(2)}:=q_2-s$ & $\alpha_6=-\alpha_1$   \\ \hline
2 & $f_7^{(1)}:=p_1, \ f_7^{(2)}:=p_2$ & $\alpha_1=0$   \\ \hline
\end{tabular}
\end{center}
Setting $\alpha_3=-\alpha_1$, then we see that the system \eqref{uraS} admits a particular solution expressed in terms of Appell's hypergeometric function. This case corresponds to the slant part of Type B in Figure 6.

\section{B{\"a}cklund transformations of the system \eqref{1}}
In this section, we study the symmetry of the Garnier system in two variables. The transformations $\pi_i \ (i=2,3,\dots,6)$ have been already obtained by H. Kimura (see \cite{K}), and the transformations $w_i,\pi_1$ have been already obtained by T. Tsuda (see \cite{Tsuda1,Tsuda2,Tsuda3,Tsuda4}.)
\begin{thm}\label{th:in1}\rm{(see \cite{K,Tsuda1,Tsuda2,Tsuda3,Tsuda4})\rm}
The system \eqref{1} admits the following transformations as its B{\"a}cklund transformations\rm{:\rm} with the notation $(*)=(q_1,p_1,q_2,p_2,t,\\
s;\alpha_1,\alpha_2,\dots,\alpha_6),$
\begin{align}\label{SG}
\begin{split}
        w_1: (*) \rightarrow& (q_1,p_1,q_2,p_2,t,s;\alpha_1+\alpha_2,-\alpha_2,\alpha_3,\alpha_4,\alpha_5,\alpha_6),\\
        w_2: (*) \rightarrow& (q_1,p_1-\frac{\alpha_3}{q_1},q_2,p_2,t,s;\alpha_1+\alpha_3,\alpha_2,-\alpha_3,\alpha_4,\alpha_5,\alpha_6),\\
        w_3: (*) \rightarrow& (q_1,p_1,q_2,p_2-\frac{\alpha_4}{q_2},t,s;\alpha_1+\alpha_4,\alpha_2,\alpha_3,-\alpha_4,\alpha_5,\alpha_6),\\
        w_4: (*) \rightarrow& (q_1,p_1-\frac{\alpha_5}{q_1+q_2-1},q_2,p_2-\frac{\alpha_5}{q_1+q_2-1},t,s;\\
        &\alpha_1+\alpha_5,\alpha_2,\alpha_3,\alpha_4,-\alpha_5,\alpha_6),\\
        w_5: (*) \rightarrow& (q_1,p_1-\frac{\alpha_6}{q_1+tq_2/s-t},q_2,p_2-\frac{\alpha_6t}{s(q_1+tq_2/s-t)},t,s;\\
        &\alpha_1+\alpha_6,\alpha_2,\alpha_3,\alpha_4,\alpha_5,-\alpha_6),\\
        \pi_1:(*) \rightarrow& (\frac{p_1(q_1p_1-\alpha_3)}{(q_1p_1+q_2p_2+\alpha_1)(q_1p_1+q_2p_2+\alpha_1+\alpha_2)},\\
        &-\frac{(q_1p_1+q_2p_2+\alpha_1)(q_1p_1+q_2p_2+\alpha_1+\alpha_2)}{p_1},\\
        &\frac{p_2(q_2p_2-\alpha_4)}{(q_1p_1+q_2p_2+\alpha_1)(q_1p_1+q_2p_2+\alpha_1+\alpha_2)},\\
        &-\frac{(q_1p_1+q_2p_2+\alpha_1)(q_1p_1+q_2p_2+\alpha_1+\alpha_2)}{p_2},\\
        &1/t,1/s;-\alpha_1-\alpha_2-\alpha_3-\alpha_4,\alpha_2,\alpha_3,\alpha_4,1-\alpha_6,1-\alpha_5),\\
        \pi_2:(*) \rightarrow& (1-q_1-q_2,-p_1,q_2,p_2-p_1,1-t,\frac{(t-1)s}{t-s};\alpha_1,\alpha_2,\alpha_5,\alpha_4,\alpha_3,\alpha_6),\\
        \pi_3:(*) \rightarrow& (\frac{-sq_1-tq_2+ts}{(t-1)s},-(t-1)p_1,\frac{(t-s)q_2}{(t-1)s},-\frac{(t-1)(tp_1-sp_2)}{t-s},\\
        &\frac{t}{t-1},\frac{t-s}{t-1};\alpha_1,\alpha_2,\alpha_6,\alpha_4,\alpha_5,\alpha_3),\\
        \pi_4:(*) \rightarrow& (\frac{q_1}{t},tp_1,\frac{q_2}{s},sp_2,\frac{1}{t},\frac{1}{s},;\alpha_1,\alpha_2,\alpha_3,\alpha_4,\alpha_6,\alpha_5),\\
        \pi_5:(*) \rightarrow& (-\frac{q_1}{q_2},-p_1q_2,\frac{1}{q_2}-(q_2p_2+q_1p_1+\alpha_1)q_2,\frac{t}{s},\frac{1}{s};\alpha_1,\alpha_4,\alpha_3,\alpha_2,\alpha_5,\alpha_6),\\
        \pi_6:(*) \rightarrow& (q_2,p_2,q_1,p_1,s,t;\alpha_1,\alpha_2,\alpha_4,\alpha_3,\alpha_5,\alpha_6).
        \end{split}
\end{align}
\end{thm}
The list \eqref{SG} should be read as
\begin{align*}
&w_1(q_1)=q_1, \quad w_1(p_1)=p_1, \quad w_1(q_2)=q_2, \quad w_1(p_2)=p_2, \quad w_1(t)=t,\\
&w_1(\alpha_1)=\alpha_1+\alpha_2, \quad w_1(\alpha_2)=-\alpha_2, \quad w_1(\alpha_3)=\alpha_3,\\
&w_1(\alpha_4)=\alpha_4, \quad w_1(\alpha_5)=\alpha_5, \quad w_1(\alpha_6)=\alpha_6.
\end{align*}

\begin{Lemma}
The transformations described in Theorem \ref{th:in1} satisfy the following relations\rm{:\rm}
\begin{align*}
&{w_1}^2={w_2}^2={w_3}^2={w_4}^2={w_5}^2=1, \quad {{\pi}_1}^2={{\pi}_2}^2={{\pi}_3}^2={{\pi}_4}^2={{\pi}_5}^2={{\pi}_6}^2=1,\\
&(w_1w_2)^2=(w_1w_3)^2=(w_1w_4)^2=(w_1w_5)^2=(w_2w_3)^2=(w_2w_4)^2\\
&=(w_2w_5)^2=(w_3w_4)^2=(w_3w_5)^2=(w_4w_5)^2=1,\\
&{\pi_1}(w_1,w_2,w_3)=(w_1,w_2,w_3){\pi_1}, \quad (w_4\pi_1)^4=(w_5\pi_1)^4=1,\\
&{\pi}_2(w_1,w_2,w_3,w_4,w_5)=(w_1,w_4,w_3,w_2,w_5){\pi}_2,\\
&{\pi}_3(w_1,w_2,w_3,w_4,w_5)=(w_1,w_5,w_3,w_4,w_2){\pi}_3,\\
&{\pi}_4(w_1,w_2,w_3,w_4,w_5)=(w_1,w_2,w_3,w_5,w_4){\pi}_4,\\
&{\pi}_5(w_1,w_2,w_3,w_4,w_5)=(w_3,w_2,w_1,w_4,w_5){\pi}_5,\\
&{\pi}_6(w_1,w_2,w_3,w_4,w_5)=(w_1,w_3,w_2,w_4,w_5){\pi}_6.
\end{align*}
\end{Lemma}

In \cite{S}, T. Suzuki showed that the system \eqref{1} has affine Weyl group symmetry of type $B_5^{(1)}$, whose generators $s_i \ (i=0,1,\dots,5)$ are explicitly given as follows: with the notation $(*)=(q_1,p_1,q_2,p_2,t,s;\gamma_0,\gamma_1,\gamma_2,\gamma_3,\gamma_4,\gamma_5),$
\begin{equation}\label{T.Suzuki}
  \left\{
  \begin{aligned}
   s_0 (*) \rightarrow &(\frac{4p_1(q_1p_1-\gamma_4-\gamma_5)}{(2q_1p_1+2q_2p_2+\gamma_0-\gamma_3-2\gamma_4-3\gamma_5)(2q_1p_1+2q_2p_2+\gamma_0+\gamma_3-\gamma_5)},\\
        &-\frac{(2q_1p_1+2q_2p_2+\gamma_0-\gamma_3-2\gamma_4-3\gamma_5)(2q_1p_1+2q_2p_2+\gamma_0+\gamma_3-\gamma_5)}{4p_1},\\
        &\frac{4p_2(q_2p_2-\gamma_5)}{(2q_1p_1+2q_2p_2+\gamma_0-\gamma_3-2\gamma_4-3\gamma_5)(2q_1p_1+2q_2p_2+\gamma_0+\gamma_3-\gamma_5)},\\
        &-\frac{(2q_1p_1+2q_2p_2+\gamma_0-\gamma_3-2\gamma_4-3\gamma_5)(2q_1p_1+2q_2p_2+\gamma_0+\gamma_3-\gamma_5)}{4p_2},\\
        &\frac{1}{t},\frac{1}{s};-\gamma_0,\gamma_1,\gamma_2+\gamma_0,\gamma_3,\gamma_4,\gamma_5),\\
        s_1 (*) \rightarrow &(\frac{q_1}{t},tp_1,\frac{q_2}{s},sp_2,\frac{1}{t},\frac{1}{s};\gamma_0,-\gamma_1,\gamma_2+\gamma_1,\gamma_3,\gamma_4,\gamma_5),\\
s_2 (*) \rightarrow &(\frac{q_1}{q_1+q_2-1},\left(p_1-(q_1p_1+q_2p_2+\frac{\gamma_0-\gamma_3-2\gamma_4-3\gamma_5}{2})\right)(q_1+q_2-1),\\
&\frac{q_2}{q_1+q_2-1},\left(p_2-(q_1p_1+q_2p_2+\frac{\gamma_0-\gamma_3-2\gamma_4-3\gamma_5}{2})\right)(q_1+q_2-1),\\
&\frac{t}{t-1},\frac{s}{s-1};\gamma_0+\gamma_2,\gamma_1+\gamma_2,-\gamma_2,\gamma_3+\gamma_2,\gamma_4,\gamma_5),\\
s_3 (*) \rightarrow &(\frac{1}{q_1},-(q_1p_1+q_2p_2+\frac{\gamma_0-\gamma_3-2\gamma_4-3\gamma_5}{2})q_1,-\frac{q_2}{q_1},-q_1p_2,\frac{1}{t},\frac{s}{t};\\
&\gamma_0,\gamma_1,\gamma_2+\gamma_3,-\gamma_3,\gamma_4+\gamma_3,\gamma_5),\\
s_4 (*) \rightarrow &(q_2,p_2,q_1,p_1,s,t;\gamma_0,\gamma_1,\gamma_2,\gamma_3+\gamma_4,-\gamma_4,\gamma_5+\gamma_4),\\
s_5 (*) \rightarrow &(q_1,p_1,q_2,p_2-\frac{\gamma_5}{q_2},t,s;\gamma_0,\gamma_1,\gamma_2,\gamma_3,\gamma_4+2\gamma_5,-\gamma_5).
   \end{aligned}
  \right. 
\end{equation}

Here, we give an explicit relation between root parameters $\gamma_0,\gamma_1,..,\gamma_5$ of type $B_5^{(1)}$ and the parameters $\alpha_1,\alpha_2,..,\alpha_6$ of the system \eqref{1} as follows:
\begin{equation}
(\gamma_0,\gamma_1,\dots,\gamma_5)=(1-\alpha_5-\alpha_6,-\alpha_5+\alpha_6,-\alpha_2+\alpha_5,\alpha_2-\alpha_3,\alpha_3-\alpha_4,\alpha_4).
\end{equation}

\begin{Lemma}
The transformations described in \eqref{T.Suzuki} define a representation of the affine Weyl group of type $B_5^{(1)}$, that is, they satisfy the following relations\rm{:\rm}
\begin{align*}
&{s_0}^2={s_1}^2={s_2}^2={s_3}^2={s_4}^2={s_5}^2=1,\\
&(s_0s_1)^2=(s_0s_3)^2=(s_0s_4)^2=(s_0s_5)^2=(s_1s_3)^2=(s_1s_4)^2=(s_1s_5)^2=\\
&(s_2s_4)^2=(s_2s_5)^2=(s_3s_5)^2=1,\\
&(s_0s_2)^3=(s_1s_2)^3=(s_2s_3)^3=(s_3s_4)^3=1, \quad (s_4s_5)^4=1.
\end{align*}
\end{Lemma}

\begin{prop}\label{B5}
The affine Weyl group $<s_0,s_1,\dots,s_5>$ of type $B_5^{(1)}$ given in \eqref{T.Suzuki} is equivalent to the group $<w_1,w_2,\dots,w_5,\pi_1,\pi_2,\dots,\pi_6>$ given in \eqref{SG}.
\end{prop}

{\bf Proof of Proposition \ref{B5}.}
We only have to check the correspondence between $s_0,s_1,\dots,s_5$ and $w_1,w_2,\dots,w_5,\pi_1,\pi_2,\dots,\pi_6$.

At first, we will show that
$$
w_1,w_2,\dots,w_5,\pi_1,\pi_2,\dots,\pi_6 \in W(B_5^{(1)}).
$$
The parameter's relations between $\alpha_i \ (i=1,2,\dots,6)$ and $\gamma_j \ (j=0,1,\dots,5)$ are given by
\begin{equation}
  \left\{
  \begin{aligned}
   \alpha_1 &=\frac{1-\gamma_1-2\gamma_2-3\gamma_3-4\gamma_4-5\gamma_5}{2},\\
   \alpha_2 &=\gamma_3+\gamma_4+\gamma_5,\\
   \alpha_3 &=\gamma_4+\gamma_5,\\
   \alpha_4 &=\gamma_5,\\
   \alpha_5 &=\gamma_1+\gamma_2+\gamma_3+\gamma_4+\gamma_5,\\
   \alpha_6 &=\gamma_2+\gamma_3+\gamma_4+\gamma_5.
   \end{aligned}
  \right. 
\end{equation}

For each $w_1,w_2,\dots,w_5,\pi_1,\pi_2,\dots,\pi_6$, we obtain the following relations:
\begin{equation}
  \left\{
  \begin{aligned}
   w_1 &=s_3s_4s_5s_4s_3,\\
   w_2 &=s_4s_5s_4,\\
   w_3 &=s_5,\\
   w_4 &=s_3s_2s_3s_4s_5s_4s_3s_2s_3,\\
   w_5 &=s_1s_3s_2s_3s_4s_5s_4s_3s_2s_3s_1,\\
   \pi_1 &=s_0,\\
   \pi_2 &=s_4s_3s_4s_3s_4s_2s_3,\\
   \pi_3 &=s_4s_3s_4s_3s_4s_2s_3s_1s_4s_3s_4s_3s_4s_2s_3,\\
   \pi_4 &=s_1,\\
   \pi_5 &=s_4s_3s_4,\\
   \pi_6 &=s_4.
   \end{aligned}
  \right. 
\end{equation}
Next, we will show that
$$
s_0,s_1,\dots,s_5 \in <w_1,w_2,\dots,w_5,\pi_1,\pi_2,\dots,\pi_6>.
$$
The parameter's relations between $\gamma_j \ (j=0,1,\dots,5)$ and $\alpha_i \ (i=1,2,\dots,6)$ are given by
\begin{equation}
  \left\{
  \begin{aligned}
   \gamma_0 &=1-\alpha_5-\alpha_6,\\
   \gamma_1 &=-\alpha_5+\alpha_6,\\
   \gamma_2 &=-\alpha_2+\alpha_5,\\
   \gamma_3 &=\alpha_2-\alpha_3,\\
   \gamma_4 &=\alpha_3-\alpha_4,\\
   \gamma_5 &=\alpha_4.
   \end{aligned}
  \right. 
\end{equation}
For each $w_1,w_2,\dots,w_5,\pi_1,\pi_2,\dots,\pi_6$, we obtain the following relations:
\begin{equation}
  \left\{
  \begin{aligned}
   s_0 &=\pi_1,\\
   s_1 &=\pi_4,\\
   s_2 &=\pi_5\pi_6\pi_5\pi_2\pi_6\pi_5\pi_6,\\
   s_3 &=\pi_6\pi_5\pi_6,\\
   s_4 &=\pi_6,\\
   s_5 &=w_3.
   \end{aligned}
  \right. 
\end{equation}
The proof has thus been completed. \qed

\section{B{\"a}cklund transformations of the system \eqref{uraS}}

\begin{thm}\label{th:in3}
The system \eqref{uraS} admits the following transformations as its B{\"a}ckl-\\
und transformations\rm{:\rm} with the notation $(*)=(q_1,p_1,q_2,p_2,t,s;\alpha_1,\alpha_2,\dots,\alpha_6),$
\begin{align}
\begin{split}
        u_1: (*) \rightarrow &(q_1+\frac{\alpha_2}{p_1},p_1,q_2,p_2,t,s;\alpha_1+\alpha_2,-\alpha_2,\alpha_3,\alpha_4,\alpha_5,\alpha_6), \\
        u_2: (*) \rightarrow &(q_1,p_1,q_2+\frac{\alpha_4}{p_2},p_2,t,s;\alpha_1+\alpha_4,\alpha_2,\alpha_3,-\alpha_4,\alpha_5,\alpha_6),\\
        u_3: (*) \rightarrow &(\frac{q_1(q_1p_1+q_2p_2-\alpha_1)}{(q_1p_1+q_2p_2-\alpha_1-\alpha_3)},\frac{p_1(q_1p_1+q_2p_2-\alpha_1-\alpha_3)}{(q_1p_1+q_2p_2-\alpha_1)},\\
        &\frac{q_2(q_1p_1+q_2p_2-\alpha_1)}{(q_1p_1+q_2p_2-\alpha_1-\alpha_3)},\frac{p_2(q_1p_1+q_2p_2-\alpha_1-\alpha_3)}{(q_1p_1+q_2p_2-\alpha_1)},\\
        &t,s;\alpha_1+\alpha_3,\alpha_2,-\alpha_3,\alpha_4,\alpha_5,\alpha_6),\\
        u_4: (*) \rightarrow &(\frac{(q_1-1)q_1p_1+(q_2-1)q_1p_2-\alpha_1q_1-\alpha_5}{(q_1-1)p_1+(q_2-1)p_2-\alpha_1-\alpha_5},\\
        &\frac{\{(q_1-1)p_1+(q_2-1)p_2-\alpha_1-\alpha_5\}p_1}{(q_1-1)p_1+(q_2-1)p_2-\alpha_1},\\
        &\frac{(q_1-1)p_1q_2+(q_2-1)q_2p_2-\alpha_1q_2-\alpha_5}{(q_1-1)p_1+(q_2-1)p_2-\alpha_1-\alpha_5},\\
        &\frac{\{(q_1-1)p_1+(q_2-1)p_2-\alpha_1-\alpha_5\}p_2}{(q_1-1)p_1+(q_2-1)p_2-\alpha_1},\\
        &t,s;\alpha_1+\alpha_5,\alpha_2,\alpha_3,\alpha_4,-\alpha_5,\alpha_6),
        \end{split}
        \end{align}
        \begin{align*}
        u_5: (*) \rightarrow &(\frac{(q_1-t)q_1p_1+(q_2-s)q_1p_2-\alpha_1q_1-\alpha_6t}{(q_1-t)p_1+(q_2-s)p_2-\alpha_1-\alpha_6},\\
        &\frac{\{(q_1-t)p_1+(q_2-s)p_2-\alpha_1-\alpha_6\}p_1}{(q_1-t)p_1+(q_2-s)p_2-\alpha_1},\\
        &\frac{(q_1-t)p_1q_2+(q_2-s)q_2p_2-\alpha_1q_2-\alpha_6s}{(q_1-t)p_1+(q_2-s)p_2-\alpha_1-\alpha_6},\\
        &\frac{\{(q_1-t)p_1+(q_2-s)p_2-\alpha_1-\alpha_6\}p_2}{(q_1-t)p_1+(q_2-s)p_2-\alpha_1},\\
        &t,s;\alpha_1+\alpha_6,\alpha_2,\alpha_3,\alpha_4,\alpha_5,-\alpha_6),\\
        \varphi_1: (*) \rightarrow &(\frac{1}{q_1},-(q_1p_1+\alpha_2)q_1,\frac{1}{q_2},-(q_2p_2+\alpha_4)q_2,\frac{1}{t},\frac{1}{s}; \\
        &-\alpha_1-\alpha_2-\alpha_3-\alpha_4,\alpha_2,\alpha_3,\alpha_4,1-\alpha_6,1-\alpha_5),\\
        \varphi_2: (*) \rightarrow &(1-q_1,-p_1,1-q_2,-p_2,1-t,1-s;\alpha_1,\alpha_2,\alpha_5,\alpha_4,\alpha_3,\alpha_6),\\
        \varphi_3: (*) \rightarrow &(\frac{t-q_1}{t-1},-(t-1)p_1,\frac{s-q_2}{s-1},-(s-1)p_2,\\
        &\frac{t}{t-1},\frac{s}{s-1};\alpha_1,\alpha_2,\alpha_6,\alpha_4,\alpha_5,\alpha_3),\\
        \varphi_4: (*) \rightarrow &(\frac{q_1}{t},tp_1,\frac{q_2}{s},sp_2,\frac{1}{t},\frac{1}{s};\alpha_1,\alpha_2,\alpha_3,\alpha_4,\alpha_6,\alpha_5),\\
        \varphi_5: (*) \rightarrow &(q_2,p_2,q_1,p_1,s,t;\alpha_1,\alpha_4,\alpha_3,\alpha_2,\alpha_5,\alpha_6),\\
        \varphi_6: (*) \rightarrow &(\frac{q_1}{q_2},p_1q_2,\frac{1}{q_2},-(q_2p_2+q_1p_1-\alpha_1)q_2,\frac{t}{s},\frac{1}{s};\alpha_1,\alpha_2,\alpha_4,\alpha_3,\alpha_5,\alpha_6),\\
        \varphi_7: (*) \rightarrow &(\frac{t}{q_1},-\frac{(q_1p_1+\alpha_2)q_1}{t},\frac{s}{q_2},-\frac{(q_2p_2+\alpha_4)q_2}{s},t,s;\\
        &-\alpha_1-\alpha_2-\alpha_3-\alpha_4,\alpha_2,\alpha_3,\alpha_4,1-\alpha_5,1-\alpha_6),\\
        \varphi_8: (*) \rightarrow &(\frac{(q_1-1)t}{q_1-t},-\frac{p_1(q_1-t)^2+\alpha_2(q_1-t)}{t(t-1)},\\
        &\frac{(q_2-1)s}{q_2-s},-\frac{p_2(q_2-s)^2+\alpha_4(q_2-s)}{s(s-1)},\\
        &t,s;\alpha_1+\alpha_3+\alpha_5-1,\alpha_2,1-\alpha_3,\alpha_4,1-\alpha_5,\alpha_6).
        \end{align*}
\end{thm}

\begin{rem}
The transformations $u_4,u_5,\varphi_7$ and $\varphi_8$ satisfy the following relations{\rm : \rm}.
\begin{align}
\begin{split}
&u_4=\varphi_2 \circ s_3 \circ \varphi_2, \quad u_5=\varphi_3 \circ s_3 \circ \varphi_3,\\
&\varphi_7=\varphi_1 \circ \varphi_4, \quad \varphi_8=\varphi_3 \circ \varphi_7 \circ \varphi_3.
\end{split}
\end{align}
\end{rem}

We note that by the transformations $S$ and $R$ (see Section 4) the relations between the transformations in Theorem \ref{th:in1} and the transformations in Theorem \ref{th:in3} are given as follows:
\begin{align}
\begin{split}
&S \circ w_i \circ R=u_i \quad (i=1,2,\dots,5),\\
&S \circ \pi_i \circ R={\varphi}_i \quad (i=1,2,\dots,6).
\end{split}
\end{align}

\begin{prop}
Let us define the following translation operators\rm{:\rm}
\begin{align}
\begin{split}
T_1:&={\varphi}_4s_4{\varphi}_1s_4,\\
T_2:&=s_4s_5T_1{\varphi}_1s_4,\\
T_3:&={\varphi}_2T_1{\varphi}_2,\\
T_4:&={\varphi}_3T_1{\varphi}_3,\\
T_5:&={\varphi}_3s_4s_5T_1{\varphi}_1s_4{\varphi}_3,\\
T_6:&={\varphi}_4s_4s_5T_1{\varphi}_1s_4{\varphi}_4,\\
T_7:&={\varphi}_1s_4s_5{\varphi}_4,\\
T_8:&={\varphi}_2s_4s_5{\varphi}_4{\varphi}_1{\varphi}_2,\\
T_9:&={\varphi}_3s_4s_5{\varphi}_4{\varphi}_1{\varphi}_3.
\end{split}
\end{align}
These translation operators act on parameters as follows\rm{:\rm}
\begin{align}
\begin{split}
T_1(\alpha_1,\alpha_2,\dots,\alpha_6)=&(\alpha_1,\alpha_2,\dots,\alpha_6)+(0,0,0,0,-1,1),\\
T_2(\alpha_1,\alpha_2,\dots,\alpha_6)=&(\alpha_1,\alpha_2,\dots,\alpha_6)+(-1,0,0,0,0,2),\\
T_3(\alpha_1,\alpha_2,\dots,\alpha_6)=&(\alpha_1,\alpha_2,\dots,\alpha_6)+(0,0,-1,0,0,1),\\
T_4(\alpha_1,\alpha_2,\dots,\alpha_6)=&(\alpha_1,\alpha_2,\dots,\alpha_6)+(0,0,1,0,-1,0),\\
T_5(\alpha_1,\alpha_2,\dots,\alpha_6)=&(\alpha_1,\alpha_2,\dots,\alpha_6)+(-1,0,2,0,0,0),\\
T_6(\alpha_1,\alpha_2,\dots,\alpha_6)=&(\alpha_1,\alpha_2,\dots,\alpha_6)+(-1,0,0,0,2,0),\\
T_7(\alpha_1,\alpha_2,\dots,\alpha_6)=&(\alpha_1,\alpha_2,\dots,\alpha_6)+(1,0,0,0,-1,-1),\\
T_8(\alpha_1,\alpha_2,\dots,\alpha_6)=&(\alpha_1,\alpha_2,\dots,\alpha_6)+(-1,0,1,0,0,1),\\
T_9(\alpha_1,\alpha_2,\dots,\alpha_6)=&(\alpha_1,\alpha_2,\dots,\alpha_6)+(-1,0,1,0,1,0).
\end{split}
\end{align}
\end{prop}

\section{Algebraic solutions of the system \eqref{uraS}}
It is known that one can get an algebraic solution of Painlev\'e VI equations by considering the fixed points with respect to a B{\"a}cklund transformation corresponding to a Dynkin automorphism.

For example, consider the Dynkin diagram automorphism
\begin{equation}
\pi_{13}(q,p,t;\alpha_0,\alpha_1,\dots,\alpha_4) \rightarrow (\frac{t}{q},-\frac{(qp+\alpha_2)q}{t},t;\alpha_0,\alpha_3,\alpha_2,\alpha_1,\alpha_4).
\end{equation}
By this transformation, the fixed solution is derived from
\begin{equation}
\alpha_0=\alpha_3, \quad \alpha_1=\alpha_4, \quad q=\frac{t}{q}, \quad p=-\frac{(qp+\alpha_2)q}{t}.
\end{equation}
Then we obtain
\begin{equation}
(q,p)=(\pm\sqrt{t},\mp\frac{\alpha_2}{2\sqrt{t}}).
\end{equation}
Masuda showed that by applying the B{\"a}cklund transformation (see \cite{TM}) defined by
\begin{equation}
T:=\pi_{13}s_3s_2s_1,
\end{equation}
we can obtain Umemura's solution with special parameters.

In this section, we consider an extension of these algebraic solutions for the system \eqref{uraS}.

At first, we consider the Dynkin diagram automorphism $\varphi_7$. By this transformation, the fixed solution is derived from
\begin{align}
\begin{split}
&\alpha_5=1-\alpha_5, \quad \alpha_6=1-\alpha_6,\\
&q_1=\frac{t}{q_1}, \quad p_1=-\frac{(q_1p_1+\alpha_2)q_1}{t}, \quad q_2=\frac{s}{q_2}, \quad p_2=-\frac{(q_2p_2+\alpha_4)q_2}{s}.
\end{split}
\end{align}
Then we obtain
\begin{align}
\begin{split}
&(\alpha_1,\alpha_2,\dots,\alpha_6)=\left(\alpha_1,\alpha_2,\alpha_3,\alpha_4,\frac{1}{2},\frac{1}{2}\right),\\
&(q_1,p_1,q_2,p_2)=\left(\pm\sqrt{t},\mp\frac{\alpha_2}{2\sqrt{t}},\pm\sqrt{s},\mp\frac{\alpha_4}{2\sqrt{s}}\right)
\end{split}
\end{align}
as a seed solution.

Applying the B{\"a}cklund transformation defined by
\begin{equation}
T:=\varphi_2\varphi_1\varphi_2u_4u_2u_1,
\end{equation}
we can obtain an extended Umemura's solution with special parameters. Here, the transformation $T$ is explicitly given as follows:
\begin{align}
\begin{split}
&T(q_1,p_1,q_2,p_2,t,s;\alpha_1,\alpha_2,\dots,\alpha_6) \rightarrow\\
&(\frac{(q_1-1)q_1p_1^2+(q_2-1)q_1p_1p_2-(\alpha_1-\alpha_2)q_1p_1}{((q_1-1)p_1+(q_2-1)p_2-\alpha_1)((q_1-1)p_1+\alpha_2)}\\
&+\frac{-\alpha_2(p_1+p_2)+\alpha_2q_2p_2-\alpha_5p_1-\alpha_1\alpha_2}{((q_1-1)p_1+(q_2-1)p_2-\alpha_1)((q_1-1)p_1+\alpha_2)},\\
&-\frac{(q_1-1)\{(q_1-1)p_1+(q_2-1)p_2-\alpha_1\}\{(q_1-1)p_1+\alpha_2\}}{(q_1-1)p_1+(q_2-1)p_2-\alpha_1-\alpha_5},\\
&\frac{q_2(q_2-1)p_2^2-(\alpha_1-\alpha_4)q_2p_2-\alpha_4(p_1+p_2)+\alpha_4q_1p_1}{((q_1-1)p_1+(q_2-1)p_2-\alpha_1)((q_2-1)p_2+\alpha_4)}\\
&+\frac{-\alpha_1\alpha_4-\alpha_5p_2+(q_1-1)p_1q_2p_2}{((q_1-1)p_1+(q_2-1)p_2-\alpha_1)((q_2-1)p_2+\alpha_4)},\\
&-\frac{(q_2-1)((q_1-1)p_1+(q_2-1)p_2-\alpha_1)((q_2-1)p_2+\alpha_4)}{(q_1-1)p_1+(q_2-1)p_2-\alpha_1-\alpha_5},\\
&\frac{t}{t-1},\frac{s}{s-1};-\alpha_1,-\alpha_2,1-\alpha_6,-\alpha_4,-\alpha_5,1-\alpha_3).
\end{split}
\end{align}

Next, we consider the Dynkin diagram automorphism $\varphi_8$. By this transformation, the fixed solution is derived from
\begin{align}
\begin{split}
&\alpha_3=1-\alpha_3, \quad \alpha_5=1-\alpha_5,\\
&q_1=\frac{(q_1-1)t}{(q_1-t)}, \quad p_1=-\frac{p_1(q_1-t)^2+\alpha_2(q_1-t)}{t(t-1)},\\
&q_2=\frac{(q_2-1)s}{(q_2-s)}, \quad p_2=-\frac{p_2(q_2-s)^2+\alpha_4(q_2-s)}{s(s-1)}.
\end{split}
\end{align}
Then we obtain
\begin{align}
\begin{split}
&(\alpha_1,\alpha_2,\dots,\alpha_6)=\left(\alpha_1,\alpha_2,\alpha_3,\alpha_4,\frac{1}{2},\frac{1}{2}\right),\\
&(q_1,p_1,q_2,p_2)=\left(t\mp\sqrt{t(t-1)},\pm\frac{\alpha_2\sqrt{t(t-1)}}{2t(t-1)},s\mp\sqrt{s(s-1)},\pm\frac{\alpha_4\sqrt{s(s-1)}}{2s(s-1)}\right)
\end{split}
\end{align}
as a seed solution.

{\it Acknowledgements.} The author would like to thank T. Masuda, T. Suzuki, N. Tahara, K. Takano, T. Tsuda and Y. Yamada  for useful discussions.

\end{document}